\documentclass[11pt]{article} 
\usepackage[francais]{babel}
\usepackage{multicol} 
\usepackage{newlfont}
\usepackage{amsmath,amssymb}

\def\zit{\mathbb{Z}}   
\def\nit{\mathbb{N}} 
\def\ppit{\mathbb{P}} 
\def\qit{\mathbb{Q}} 
\def\cit{\mathbb{C}} 

\newcommand{\pf}{{\bf Preuve.~}}

\newcommand{\qed}{\hfill~~\mbox{$\Box$}}

\newenvironment{proof}{\smallskip \noindent \pf}{\qed \bigskip}

\newtheorem{theorem}{Th\'eor\`eme}[subsection]
\newtheorem{proposition}[theorem]{Proposition}

\newtheorem{lemma}[theorem]{Lemme}
\newtheorem{corollary}[theorem]{Corollaire}

\newtheorem{remark}[theorem]{Remarque}
\newtheorem{example}[theorem]{Exemple}

\begin{document}

\title{{\bf Construction de vari\'et\'es de Frobenius {\em via} les polyn\^omes de Laurent : une autre approche}
\author{\sc Antoine Douai}
\thanks{Mots clef : Syst\`emes de Gauss-Manin,r\'eseau de Brieskorn, vari\'et\'es de Frobenius.
Classification AMS : 32S40, 32S30.}}

\maketitle

\section{Introduction}

Soit $f$ un polyn\^ome de Laurent non d\'eg\'en\'er\'e 
par rapport \`a son poly\`edre de Newton \`a l'infini et commode. 
Il est maintenant \'etabli (il s'agit l\`a du r\'esultat principal de [DoSa]) que l'on peut munir l'espace des param\`etres $\tilde{M}$ du d\'eploiement 
universel d'un tel objet d'une structure de Frobenius : $\tilde{M}$ est ainsi une {\em vari\'et\'e de Frobenius}. 
Le premier ingr\'edient pour obtenir une telle structure est un fibr\'e trivial ${\cal G}$ sur $\ppit^{1}\times\tilde{M}$, muni 
d'une connexion m\'eromorphe plate \`a p\^oles de type $1$ le long de $\{\theta =0\}\times\tilde{M}$ et logarithmique le long 
de $\{\theta =\infty \}\times\tilde{M}$, la carte centr\'ee \`a l'origine \'etant munie de la coordonn\'ee $\theta$. Le second est 
une forme primitive et homog\`ene qui permet 
d'identifier la restriction ${\cal G}_{|\{0\}\times\tilde{M}}$
au fibr\'e tangent $T\tilde{M}$ et de transporter ainsi sur $T\tilde{M}$ les diff\'erentes structures naturellement construites sur 
${\cal G}_{|\{0\}\times\tilde{M}}$.
Classiquement (c'est la m\'ethode propos\'ee par K. Saito et reprise dans [DoSa]), le fibr\'e ${\cal G}$
s'obtient en r\'esolvant le probl\`eme de Birkhoff pour le r\'eseau de Brieskorn 
du d\'eploiement universel du polyn\^ome consid\'er\'e, fibr\'e intrins\`equement muni d'une connexion m\'eromorphe $\nabla$ de type $1$ le long de
$\{\theta =0\}\times\tilde{M}$. Ceci peut \^etre fait pour n'importe quel polyn\^ome de Laurent 
$M$-mod\'er\'e, c'est \`a dire v\'erifiant une propri\'et\'e 
de mod\'eration \`a l'infini topologique, et pas seulement pour les polyn\^omes commodes et non d\'eg\'en\'er\'es. 
Par contre, l'existence d'une forme primitive {\em et} homog\`ene n'est pour le moment garantie que dans 
le cas non d\'eg\'en\'er\'e [DoSa, 4.d], ce qui explique 
l'hypoth\`ese plus restrictive faite sur $f$. Contrairement au cas local, c'est \`a dire lorsque $f$ est un germe de fonction 
holomorphe \`a point critique isol\'e \`a l'origine, un \'enonc\'e plus g\'en\'eral semble ici sans espoir.

La premi\`ere \'etape, \`a savoir l'analyse du r\'eseau de Brieskorn du d\'eploiement universel du polyn\^ome $f$, est une difficult\'e 
majeure : en effet, sa coh\'erence, finalement garantie par l'hypoth\`ese de mod\'eration faite sur $f$, se r\'ev\`ele particuli\`erement d\'elicate. 
De plus, les inopportuns points critiques qui disparaissent \`a l'infini imposent dans un premier temps  
le recours \`a l'analytisation dans les variables du polyn\^ome et par la suite l'emploi d'un proc\'ed\'e d'alg\'ebrisation transcendant.  
Ces points sont d\'etaill\'es dans les deux premiers paragraphes de [DoSa] et font appel \`a des techniques et des r\'esultats plut\^ot 
sophistiqu\'es.  
L'objet de cet article est de montrer 
comment on peut les \'eviter 
en utilisant un joli r\'esultat de C. Hertling et Y. Manin [HeMa, theorem 4.5] qui, nous le verrons, s'adapte miraculeusement au cas polynomial : il n'est en fait pas n\'ecessaire de consid\'erer le 
d\'eploiement universel du polyn\^ome $f$, comme cela est fait dans [DoSa], pour obtenir 
le fibr\'e trivial ${\cal G}$ attendu. Dans le cadre alg\'ebrique, on obtient ainsi une m\'ethode de construction de vari\'et\'es de Frobenius dont la 
pr\'esentation est plus \'el\'ementaire et autonome.
 
Le point de d\'epart est l'\'etude du r\'eseau de Brieskorn $G_{0}$ 
 associ\'e \`a la d\'eformation $F$ d'un polyn\^ome de Laurent $f$ non d\'eg\'en\'er\'e et commode d\'efinie par
$$F(u,x)=f(u)+\sum_{j=1}^{r}x_{j}g_{j}(u)$$
o\`u $u=(u_{1},\cdots ,u_{n})$ sont les coordonn\'ees sur $U=(\cit^{*})^{n}$, 
$x=(x_{1},\cdots ,x_{r})\in\cit^{r}$ les param\`etres de la d\'eformation et les $g_{j}$ 
sont des polyn\^omes de Laurent en $u$. Nous supposons que les polyn\^omes $g_{j}$ sont 
''strictement sous diagramme'', c'est \`a dire s'\'ecrivent $\sum a_{\alpha}u^{\alpha}$ avec $\alpha$ strictement \`a l'int\'erieur du poly\`edre de Newton de $f$. $F$ d\'efinit alors une famille de polyn\^omes de Laurent non d\'eg\'en\'er\'es et commodes \`a nombre de Milnor {\em global} constant. 
Nous montrons alors, en utilisant des m\'ethodes \'el\'ementaires, 
que $G_{0}$ est libre, de type fini sur $\cit [x,\theta ]$. 
Il s'agit l\`a du r\'esultat principal du paragraphe 2
et il faut de se r\'ef\'erer \`a [DoSa, 2.b] pour mesurer le gain obtenu. 
La r\'esolution du probl\`eme de Birkhoff pour 
$G_{0}$, abord\'ee au paragraphe 3, 
 peut alors se traiter alg\'ebriquement. Pr\'ecis\'ement,
nous montrons qu'il existe une base de $G_{0}$ dans laquelle la matrice de la connexion plate $\nabla$ s'\'ecrit
$$(\frac{B_{0}(x)}{\theta}+B_{\infty})\frac{d\theta }{\theta}+\sum_{i=1}^{r}\frac{C^{(i)}(x)}{\theta}dx_{i}$$
\noindent{o\`u} $B_{\infty}$ est une matrice constante, les matrices $B_{0}(x)$ et $C^{(i)}(x)$ \'etant \`a coefficients dans $\cit [x]$.
Bien s\^ur, par unicit\'e, la solution du probl\`eme de Birkhoff obtenue ici coincide avec celle donn\'ee 
par le th\'eor\`eme g\'en\'eral de B. Malgrange [Mal, th\'eor\`eme 2.2], mais notre m\'ethode 
\`a l'avantage d'\^etre plus explicite (voir l'exemple 4.2.2). En particulier,
la restriction de cette solution \`a toute valeur du param\`etre est \'egale 
\`a la solution de la restriction donn\'ee par [DoSa, appendix B]. 
Nous remarquons aussi que l'on peut d\'efinir, sur le syst\`eme de Gauss-Manin d'une telle d\'eformation, 
une forme sesquilin\'eaire 
non d\'eg\'en\'er\'ee, de poids $n$ relativement au r\'eseau $G_{0}$ : cette dualit\'e joue un r\^ole essentiel dans la construction des vari\'et\'es 
de Frobenius.
En d'autres termes, nous construisons une 
$trTLEP(n)$-structure au sens de [HeMa, definition 4.1] qui permet \`a son tour (voir [HeMa, theorem 4.2] et les r\'ef\'erences 
qui y sont mentionn\'ees) de d\'efinir une structure de type Frobenius, c'est \`a dire un 
uplet
$$(M, E, \bigtriangledown ,\Phi ,R_{0}, R_{\infty},g)$$
o\`u $E$ est un fibr\'e sur $M$, $\bigtriangledown$ une connexion plate sur $E$, $R_{0}$ et $R_{\infty}$ des endormorphismes de $E$, $\Phi$ un champ de Higgs et $g$ une m\'etrique sur $E$. Tous ces objets v\'erifient des relations de compatibilit\'e naturelles. Ce sont les ``conditions initiales'' qui vont permettre d'obtenir la structure de Frobenius sur $\tilde{M}$.

Nous remarquons ensuite, en utilisant [HeMa, theorem 4.5], que cette structure de type Frobenius se d\'eploie de 
mani\`ere universelle (de fa\c{c}on toutefois transcendante, 
mais cette fois ci le proc\'ed\'e, auquel il semble difficile d'\'echapper, concerne {\em seulement les param\`etres de la d\'eformation}), pourvu que les polyn\^omes $g_{i}$ soient choisis de sorte que
$$\cit [g_{1},\cdots ,g_{r}].\zeta =
\frac{\Omega^{n}(U)}{df\wedge \Omega^{n-1}(U)}:=E_{0},$$
$\zeta$ d\'esignant la classe de la forme $\frac{du_{1}}{u_{1}}\wedge\cdots \wedge \frac{du_{n}}{u_{n}}$ dans $E_{0}$ et $\cit [g_{1},\cdots ,g_{r}]$ la sous alg\`ebre des endomorphismes de $E_{0}$ engendr\'ee par la multiplication par les $g_{i}$. Ceci fournit la condition (GC) de {\em loc. cit.}, la condition (IC) \'etant alors essentiellement garantie par le fait que les $g_{i}$ sont strictement sous diagramme. 
L'existence d'une structure de Frobenius sur $\tilde{M}$  
est alors assujettie \`a l'homog\'en\'eit\'e de la forme $\zeta$ : nous retrouvons ici la situation classique. 
Ceci est expliqu\'e au paragraphe 4. 
Signalons que la structure de Frobenius obtenue ici d\'epend donc au premier abord de deux choix, celui de la d\'eformation  $F$ et celui de la 
forme $\zeta$, qu'il serait agr\'eable de fixer une fois pour toute. Ce probl\`eme reste ouvert. Plus g\'en\'eralement, et c'est l\`a une particularit\'e inattendue du cas alg\'ebrique, une forme primitive et homog\`ene n'est pas n\'ecessairement unique. 

Si l'on part maintenant d'une fonction polynomiale, commode et non d\'eg\'en\'er\'ee, c'est \`a dire si l'on travaille sur $U=\cit^{n}$, la m\'ethode employ\'ee ci-dessus reste pertinente, une
fois une forme homog\`ene $\zeta$ v\'erifiant les conditions de [HeMa, theorem 4.5] trouv\'ee, m\^eme si certains \'enonc\'es doivent \^etre l\'eg\`erement modifi\'es. Elle s'applique id\'ealement lorsque $f$ est un polyn\^ome quasi-homog\`ene, cas pour lequel on a $R_{0}=0$ sur $E_{0}$. Ceci est discut\'e 
\`a la fin du paragraphe 4. Nous terminons ce chapitre par un exemple d\'etaill\'e. Il est essentiel de noter que la situation alg\'ebrique globale consid\'er\'ee ici, et pour laquelle le r\'esultat de Manin et Hertling est si bien adapt\'e, est diff\'erente de celle, analytique et locale, traditionnellement trait\'ee : nous regardons des d\'eformations de polyn\^omes (de Laurent) \`a nombre de Milnor global constant et non pas des d\'eformations d'un point critique \`a nombre de Milnor local constant.

Ind\'ependamment de la construction de vari\'et\'es de Frobenius, il peut \^etre int\'eressant d'\'etudier  plus 
g\'en\'eralement le r\'eseau de Brieskorn d'une famille de  
polyn\^omes mod\'er\'es
\`a nombre de Milnor global constant et le probl\`eme de Birkhoff aff\'erent. L'\'etude du cas non d\'eg\'en\'er\'e, couverte par les 
paragraphes 2 et 3 de ces notes, donne une id\'ee des r\'esultats auxquels on doit s'attendre.
Le paragraphe 5 explore les pistes pour 
y parvenir.\\ 

Je remercie C. Sabbah qui m'a signal\'e le r\'esultat de C. Hertling et Y. Manin sur lequel 
repose cet article mais aussi
pour de fructueuses discussions, en particulier lors de l'\'elaboration de ce texte. Je remercie aussi C. Hertling pour sa lecture attentive d'une premi\`ere version de ce travail.\\

{\bf Notations :} Dans ce qui suit, $U=(\cit^{*})^{n}$. $U$ est muni des coordonn\'ees $u=(u_{1},\cdots ,u_{n})$ et
on note $K=\cit [u_{1},\cdots u_{n},u_{1}^{-1},\cdots ,u_{n}^{-1}]$.\\

\section{Le r\'eseau de Brieskorn : le cas sous-diagramme}
Soit  
$f:U\rightarrow\cit$ un polyn\^ome de Laurent et 
$$F:U\times\cit^{r}\rightarrow\cit$$
la d\'eformation de $f$ d\'efinie par, en notant $x=(x_{1},\cdots ,x_{r})$ les coordonn\'ees sur $\cit^{r}$,
$$F(u,x)=f(u)+\sum_{j=1}^{r}x_{j}g_{j}(u)$$
les $g_{j}$ d\'esignant des polyn\^omes de Laurent donn\'es.
Nous supposerons $f$ {\em non d\'eg\'en\'er\'e par rapport \`a son poly\`edre de Newton 
\`a l'infini et commode} [K]. En particulier, $f$ n'a que des points critiques isol\'es et on d\'efinit le {\em nombre de Milnor global} $\mu$ par la formule
$$\mu =\dim_{\cit}\frac{\Omega^{n}(U)}{df\wedge\Omega^{n-1}(U)}$$
\noindent{o\`u} $\Omega^{p}(U)$ d\'esigne l'espace des $p$-formes r\'eguli\`eres sur $U$.
On d\'efinit 
$$G_{0}=\frac{\Omega^{n}(U)[x,\theta ]}{(\theta d_{u}-d_{u}F\wedge) \Omega^{n-1}(U)[x,\theta ]}$$
(la notation $d_{u}$ signifie que la diff\'erentielle n'est prise que par rapport \`a $u$)
qui est naturellement un $\cit [x,\theta ]$-module que nous appellerons {\em module de Brieskorn} et le $\cit [x,\theta ,\theta^{-1}]$-module 
$$G=\Omega^{n}(U)[x,\theta ,\theta^{-1}]/(\theta d_{u}-d_{u}F\wedge)\Omega^{n-1}(U)[x,\theta ,\theta^{-1}]$$
le localis\'e du transform\'e de Fourier du syst\`eme de Gauss-Manin : nous l'appellerons aussi syst\`eme de Gauss-Manin.
On pose $\tau :=\theta^{-1}$ et, pour $\omega\in\Omega^{n}(U)[x]$,
$$\nabla_{\partial_{\tau}}(\omega\tau^{p}):=-F\omega\tau^{p}+p\omega\tau^{p-1},$$
et
$$\nabla_{\partial_{x_{j}}}(\omega\tau^{p}):=\partial_{x_{j}}(\omega )\tau^{p}-\frac{\partial F}{\partial x_{j}}\omega\tau^{p+1}$$
pour tout $j\in \{1,\cdots ,r\}$. Ces deux actions commutant avec la diff\'erentielle $\theta d_{u}-d_{u}F\wedge$, le syst\`eme de Gauss-Manin 
$G$ se trouve muni d'une connexion $\nabla$. Remarquons qu'alors $G_{0}$ est stable par $\theta^{2}\nabla_{\partial_{\theta}}$.
On notera aussi
$$G_{0}^{o}=\frac{\Omega^{n}(U)[\theta ]}{(\theta d-df\wedge) \Omega^{n-1}(U)[\theta ]}$$
et
$$G^{o}=\Omega^{n}(U)[\theta ,\theta^{-1}]/(\theta d-df\wedge)\Omega^{n-1}(U)[\theta ,\theta^{-1}].$$

\noindent{L'objet} de ce paragraphe est d'\'etendre en famille, dans ce cas pr\'ecis, 
les r\'esultats relatifs au r\'eseau de Brieskorn d'un polyn\^ome de Laurent $f$  non d\'eg\'en\'er\'e par rapport \`a son poly\`edre de Newton 
\`a l'infini et commode expos\'es dans [DoSa, section 4].

\subsection{Pr\'eliminaires}

\subsubsection{Filtrations}

On d\'efinit tout d'abord une filtration de Newton sur $G_{0}$ et $G$.
On garde les notations de [DoSa, 4.a] o\`u il est construit une filtration de Newton ${\cal N}_{\bullet}$ sur $K$ index\'ee par $\qit$.
On d\'efinit de mani\`ere similaire une filtration, not\'ee aussi $\cal{N}_{\bullet}$, sur $K[x]$ en donnant aux param\`etres $x_{i}$ le poids $0$ : pour 
$g\in K[x]$ et $\sigma$ une face de dimension $n-1$ de la fronti\`ere de Newton, on pose $\phi_{\sigma}(g)=max_{a\in supp(g)}L_{\sigma}(a)$ o\`u 
$L_{\sigma}$ d\'esigne la forme lin\'eaire \`a coefficients rationnels (pas n\'ecessairement positifs) qui prend la valeur $1$ sur $\sigma$, le maximum 
\'etant pris sur les
exposants des mon\^omes {\em en $u$} apparaissant dans $g$; on pose alors $\phi (g)=max_{\sigma}\phi_{\sigma}(g)$ et
$${\cal N}_{\alpha}K[x]:=\{v\in K[x], \phi (v)\leq\alpha \}.$$
\noindent{On} en d\'eduit une filtration ${\cal N}_{\bullet}$ sur $\Omega^{n}(U)[x]$ en posant
$${\cal N}_{\alpha}\Omega^{n}(U)[x]:=\{v\frac{du}{u}\in \Omega^{n}(U)[x], \phi (v)\leq\alpha \}$$
\noindent o\`u $du/u =du_{1}/u_{1}\wedge \cdots \wedge du_{n}/u_{n}$, puis sur $\Omega^{n}(U)[x,\theta ]$ en posant 
$${\cal N}_{\alpha}\Omega^{n}(U)[x,\theta ]:={\cal N}_{\alpha}\Omega^{n}(U)[x]+\theta {\cal N}_{\alpha -1}\Omega^{n}(U)[x]+\cdots +\theta^{k}{\cal N}_{\alpha -k}\Omega^{n}(U)[x]+\cdots $$ 
\noindent Remarquons que la filtration ${\cal N}_{\bullet}$ sur $\Omega^{n}(U)[x]$ est normalis\'ee de sorte que $du/u\in {\cal N}_{0}\Omega^{n}(U)[x]$.
On obtient, par projection, une filtration ${\cal N}_{\bullet}$ sur $G_{0}$ ({\em cf.} [DoSa, definition 4.2]) et
on d\'efinit une filtration ${\cal N}_{\bullet}$ sur $G$ en posant, pour tout $\alpha\in\qit$,
$${\cal N}_{\alpha}G:={\cal N}_{\alpha}G_{0}+\tau{\cal N}_{\alpha +1}G_{0}+\cdots +\tau^{k}{\cal N}_{\alpha +k}G_{0}+\cdots$$

\subsubsection{Un lemme de division}

Pour $\beta\in\qit$, soit $E_{\beta}\subset {\cal N}_{\beta}K\subset {\cal N}_{\beta}K[x]$ 
un sous espace vectoriel isomorphe  par projection \`a
$$\frac{{\cal N}_{\beta}K}{(\partial f)\cap {\cal N}_{\beta}K+{\cal N}_{<\beta}K},$$
$(\partial f)$ d\'esignant l'id\'eal engendr\'e par les d\'eriv\'ees partielles de $f$,
et $E_{\beta}[x]$ le $\cit [x]$-module libre engendr\'e par $E_{\beta}$.
Rappelons que l'on a suppos\'e $f$ commode et non d\'eg\'en\'er\'e.

\begin{lemma} On suppose que, pour tout $i\in \{1,\cdots ,n\}$ et tout $j\in \{1,\cdots ,r\}$, on a 
$$\phi (u_{i}\frac{\partial g_{j}}{\partial u_{i}})<1.$$
 Alors, 
pour tout $h\in{\cal N}_{\alpha}K[x]$, il existe $v_{\alpha}\in \oplus_{\beta\leq\alpha}E_{\beta}[x]$ et $a_{1},\cdots ,a_{n}\in K[x]$ tels que
$$h=v_{\alpha}+\sum_{i=1}^{n}a_{i}u_{i}\frac{\partial F}{\partial u_{i}}$$
les $a_{i}$ v\'erifiant de plus $\phi (a_{i})\leq \alpha -1$ et 
$\phi (\frac{\partial u_{i}a_{i}}{\partial u_{i}})\leq \alpha -1$.
\end{lemma}
\begin{proof} On utilise le th\'eor\`eme de division absolu de Kouchnirenko [K, th\'eor\`eme 4.1]
comme dans [DoSa, paragraphe 4.b] : 
si $h\in{\cal N}_{\alpha}K$, il existe $w_{\alpha}\in\oplus_{\beta\leq\alpha}E_{\beta}$ et
$b_{1},\cdots , b_{n}\in K$ tels que
$$h=w_{\alpha}+\sum_{i=1}^{n}b_{i}u_{i}\frac{\partial f}{\partial u_{i}}$$
avec de plus $\phi (b_{i})\leq \alpha -1$ et $\phi (\frac{\partial u_{i}b_{i}}{\partial u_{i}})\leq \alpha -1$. On en d\'eduit que 
$$h=w_{\alpha}+\sum_{i=1}^{n}b_{i}u_{i}\frac{\partial F}{\partial u_{i}}-\sum_{i=1}^{n}b_{i}[\sum_{j=1}^{r}x_{j}u_{i}\frac{\partial g_{j}}{\partial u_{i}}].$$
Posons $g:=\sum_{j=1}^{r}x_{j}g_{j}$. On a 
$$\phi (b_{i}u_{i}\frac{\partial g}{\partial u_{i}})\leq \phi (b_{i})+\phi (u_{i}\frac{\partial g}{\partial u_{i}})\leq 
\alpha -1+\phi ( u_{i}\frac{\partial g}{\partial u_{i}}),$$
la premi\`ere in\'egalit\'e r\'esultant des propri\'et\'es de $\phi$ et la seconde des conditions sur les poids impos\'ees \`a $b_{i}$. Par hypoth\`ese, on a
$\phi (u_{i}\frac{\partial g}{\partial u_{i}})<1$ et l'on en d\'eduit que 
$$\phi (b_{i}u_{i}\frac{\partial g}{\partial u_{i}})<\alpha.$$
Cette in\'egalit\'e stricte permet de conclure par r\'ecurrence parce que ${\cal N}_{\alpha}K=0$ pour $\alpha <0$.
\end{proof} 

\subsection{Finitude de $G_{0}$ sur $\cit [x,\theta ]$}

Rappelons que
$$F(u,x)=f(u)+\sum_{j=1}^{r}x_{j}g_{j}(u),$$
$f$ \'etant un polyn\^ome de Laurent commode et non d\'eg\'en\'er\'e et que $G_{0}$ d\'esigne le module de Brieskorn associ\'e \`a $F$. 

\begin{proposition} On suppose que l'on a $\phi (g_{j})<1$ pour tout $j\in \{1,\cdots ,r\}$. Alors $G_{0}$ est libre de type fini, 
de rang $\mu$ sur $\cit [x,\theta ]$.
\end{proposition} 
\begin{proof} En effet, l'hypoth\`ese sur les $g_{j}$ montre que $\phi (u_{i}\frac{\partial g_{j}}{\partial u_{i}})<1$ pour tout $i\in \{1,\cdots ,n\}$ et tout $j\in \{1,\cdots ,r\}$. 
Le lemme de division 2.1.1 permet alors de montrer, de la m\^eme mani\`ere qu'en [DoSa, remark 4.8] et en fournissant un syst\`eme 
de g\'en\'erateurs adapt\'e \`a la filtration de Newton, que
$G_{0}$ est de type fini sur $\cit [x,\theta ]$. On dispose donc d'un morphisme surjectif
$$\cit [x,\theta ]^{\mu}\rightarrow G_{0}.$$
Soit $N$ son noyau. Pour toute valeur fix\'ee du param\`etre $x^{0}=(x_{1}^{0},\cdots ,x_{r}^{0})$, 
$\oplus_{\alpha}E_{\alpha}$ fournit aussi une base sur $\cit [\theta]$
du r\'eseau de Brieskorn associ\'e au polyn\^ome $F(u,x^{0})=f(u)+\sum_{j=1}^{r}x_{j}^{0}g_{j}(u)$ : en effet,
la condition $\phi (g_{j})<1$ entra\^ine que, pour toute valeur $x^{0}=(x^{0}_{1},\cdots ,x^{0}_{n})$ fix\'ee du param\`etre $x$, le polyn\^ome 
de Laurent $F(u,x^{0})=f(u)+\sum_{j=1}^{r}x^{0}_{j}g_{j}$ est commode et non d\'eg\'en\'er\'e (son nombre de Milnor global \'etant
 de plus \'egal au nombre 
de Milnor global de $f$) et que, avec les notations du d\'ebut du paragraphe 2.1.2,
$$(\partial f)\cap {\cal N}_{\alpha}K+{\cal N}_{<\alpha}K=(\partial F(.,x^{0}))\cap {\cal N}_{\alpha}K+{\cal N}_{<\alpha}K.$$
 Un \'el\'ement de $N$ est donc repr\'esent\'e par un nombre fini de polyn\^omes en $\theta$ 
dont les coefficients sont des polyn\^omes de  $\cit [x]$ qui s'annulent en tout point. Ainsi $N$ est-il nul.
\end{proof}

\begin{corollary} On suppose que l'on a $\phi (g_{j})<1$ pour tout $j\in \{1,\cdots ,r\}$. 
Alors $G$ est un $\cit [x,\tau ,\tau^{-1}]$-module libre de rang $\mu$.
\end{corollary}

\noindent Sous les hypoth\`eses de la proposition 2.2.1, $G_{0}$ est un r\'eseau de $G$ : c'est le r\'eseau de Brieskorn. 

\subsection{$V$-filtration et modules de Hodge}
Nous consid\'erons maintenant la situation suivante : la d\'eformation $F$ du polyn\^ome de Laurent commode et non d\'eg\'en\'er\'e $f$ est, comme ci-dessus,
$$F(u,x)=f(u)+\sum_{j=1}^{r}x_{j}g_{j}(u)$$
et nous faisons l'hypoth\`ese suppl\'ementaire que {\em $\phi (g_{j})<1$ pour tout $j\in \{1,\cdots ,r\}$}. La proposition 2.2.1 et le corollaire 2.2.2 s'appliquent donc.\\

\noindent On d\'efinit, pour $\alpha\in\qit$ et $p\in\zit$,
$$H_{\alpha}:={\cal N}_{\alpha}G/ {\cal N}_{<\alpha}G,$$
$$F_{p}H_{\alpha}:=({\cal N}_{\alpha}G\cap\tau^{p}G_{0}+{\cal N}_{<\alpha}G)/{\cal N}_{<\alpha}G$$
et leurs homologues pour la valeur nulle des param\`etres
$$H_{\alpha}^{o}:={\cal N}_{\alpha}G^{o}/ {\cal N}_{<\alpha}G^{o},$$
$$F_{p}H_{\alpha}^{o}:=({\cal N}_{\alpha}G^{o}\cap\tau^{p}G_{0}^{o}+{\cal N}_{<\alpha}G^{o})/{\cal N}_{<\alpha}G^{o}.$$

\noindent D'apr\`es la proposition 2.2.1,
le module de Brieskorn $G_{0}$ est libre de type fini sur 
$\cit [x,\theta ]$  et
le lemme de division 2.1.1 a d'autres cons\'equences, r\'eunies dans la

\begin{proposition} On suppose que $\phi (g_{j})<1$ pour tout $j$. Alors,\\
1. ${\cal N}_{\alpha}G_{0}={\cal N}_{\alpha}G\cap G_{0}$ pour tout $\alpha\in\qit$,\\
2. $\nabla_{\partial_{x_{i}}}({\cal N}_{\alpha}G)\subset {\cal N}_{\alpha}G$ pour tout $i\in \{1,\cdots ,r\}$,\\
3. pour tout $\alpha\in\qit$, ${\cal N}_{\alpha}G_{0}$ est un $\cit[x]$-module libre de type fini,\\
4. pour tout $\alpha\in\qit$, ${\cal N}_{\alpha}G$ est un $\cit [x,\tau ]$-module de type fini,\\
5. pour tout $\alpha\in\qit$, $H_{\alpha}$ est un $\cit [x]$-module libre de type fini,\\ 
6. pour tout $p\in\zit$, $\alpha\in\qit$ et $i\in \{1,\cdots ,r\}$,
$$\nabla_{\partial_{x_{i}}}F_{p}H_{\alpha}\subset F_{p}H_{\alpha},$$\\
7. pour tout $p\in\zit$ et $\alpha\in\qit$, 
$F_{p}H_{\alpha}$ est un $\cit [x]$-module libre de type fini.
\end{proposition}
\begin{proof}
1. est montr\'e comme en [DoSa, End of the proof of theorem 4.5, p. 1093].\\
2. On a 
$$\nabla_{\partial_{x_{i}}}(\omega\tau^{p})=\partial_{x_{i}}(\omega )\tau^{p}-\frac{\partial F}{\partial x_{i}}\omega\tau^{p+1}$$ 
et il suffit de montrer que $\nabla_{\partial_{x_{i}}}{\cal N}_{\alpha}G_{0}\subset \tau ({\cal N}_{\alpha +1}G_{0})$.
Si $\omega\in{\cal N}_{\alpha}G_{0}$, on a clairement 
$$\partial_{x_{i}}(\omega )\subset {\cal N}_{\alpha}G_{0}.$$ 
Le r\'esultat s'obtient en remarquant que
$\frac{\partial F}{\partial x_{i}}\omega\tau = g_{i}\omega\tau$ pour 
$i\in \{1,\cdots ,r\}$ et en se souvenant que, par hypoth\`ese, 
$$g_{i}\omega \in {\cal N}_{<1+\alpha}G_{0}.$$
On a en particulier  
$$\frac{\partial F}{\partial x_{i}}\omega\tau\in \tau ({\cal N}_{<\alpha +1}G\cap G_{0})\subset {\cal N}_{<\alpha}G.$$
3. La finitude s'obtient en utilisant le lemme 2.1.1 comme en [DoSa, remark 4.8]; la libert\'e s'obtient comme dans la preuve de 
la proposition 2.2.1.\\
4. D'apr\`es le lemme 2.1.1, pour chaque $\alpha$ il existe $k_{0}$ tel que
$${\cal N}_{\alpha}G\subset\cit [x,\tau ]({\cal N}_{\alpha}G_{0}+\cdots +\tau^{k_{0}}{\cal N}_{\alpha +k_{0}}G_{0})$$
et chaque $\tau^{i}{\cal N}_{\alpha +i}G_{0}$ est un $\cit [x]$-module de type fini d'apr\`es 3.\\
5. En effet, 4. montre que 
$${\cal N}_{\alpha}G\subset \cit [x]({\cal N}_{\alpha}G_{0}+\cdots +
\tau^{k_{0}}{\cal N}_{\alpha +k_{0}}G_{0})+
{\cal N}_{<\alpha}G.$$\\
6. S'obtient \`a l'aide du point 2. qui montre que 
$$\frac{\partial F}{\partial x_{i}}\omega\tau\subset {\cal N}_{<\alpha}G.$$\\
7. $F_{p}H_{\alpha}$ est de type fini car $H_{\alpha}$ l'est; de 6. on d\'eduit alors qu'il est libre. 
\end{proof}

\begin{remark} En g\'en\'eral, on a seulement la condition de transversalit\'e
$$\nabla_{\partial_{x_{i}}}F_{p}H_{\alpha}\subset F_{p+1}H_{\alpha}.$$
L'hypoth\`ese $\phi (g_{i})<1$ est essentielle pour obtenir le point 6.
\end{remark}

\noindent{Nous} pouvons maintenant \'enoncer

\begin{proposition} Sous les hypoth\`eses du d\'ebut du paragraphe 2.3, la filtration
${\cal N}_{\bullet}G$ est \'egale \`a la $V$-filtration de Malgrange-Kashiwara de $G$ le long de $\tau =0$.
\end{proposition}
\begin{proof}
Posons $D_{X}:=\cit [x_{1},\cdots ,x_{r}]<\partial_{x_{1}},\cdots ,\partial_{x_{r}}>$. On
d\'efinit une filtration 
$V_{\bullet}$ sur $D_{X}[\tau ]<\partial_{\tau}>$, index\'ee par $\zit$, en posant
$$V_{0}D_{X}[\tau ]<\partial_{\tau}>=D_{X}[\tau ]<\tau\partial_{\tau}>,$$
$$V_{k}D_{X}[\tau ]<\partial_{\tau}>=V_{k-1}+\partial_{\tau}V_{k-1}\ pour\ k\geq 1,$$
$$V_{k}D_{X}[\tau ]<\partial_{\tau}>=\tau^{-k}V_{0}D_{X}[\tau ]<\partial_{\tau}>\ pour\ k\leq 0.$$
Par unicit\'e, il suffit de montrer que \\
i) chaque ${\cal N}_{\alpha}G$ est un $V_{0}D_{X}[\tau ]<\partial_{\tau}>$-module de type fini,\\
ii) $V_{i}D_{X}[\tau ]<\partial_{\tau}>({\cal N}_{\alpha}G)\subset {\cal N}_{\alpha +i}G$ pour tout $\alpha\in\qit$ et tout $i\in\zit$,\\
iii) pour tout $\alpha\in\qit$, l'action de $\tau\partial_{\tau}+\alpha$ est nilpotente sur 
$Gr^{{\cal N}}_{\alpha}G:={\cal N}_{\alpha}G/ {\cal N}_{<\alpha}G.$\\
Observons tout d'abord que ${\cal N}_{\bullet}$ est stable par $\tau\partial_{\tau}$. Ceci r\'esulte de la d\'efinition de ${\cal N}_{\alpha}G$ donn\'ee en 2.1.1 : si $\omega\tau^{i}\in {\cal N}_{\alpha}G$, on a en effet
$$\tau\partial_{\tau}(\omega\tau^{i})=-F\omega\tau^{i+1}+i\omega\tau^{i}$$
avec $F\omega\tau^{i+1}\in {\cal N}_{\alpha}G$ parce que $F\omega\in {\cal N}_{\alpha +1}G_{0}$ et $\omega\tau^{i}\in {\cal N}_{\alpha}G.$
Le point i) r\'esulte alors du point 4. de la proposition 2.3.1 et les points ii) et iii) se montrent maintenant de la m\^eme mani\`ere
qu'en [DoSa, lemma 4.11], l'outil essentiel \'etant le lemme de division 2.1.1.
\end{proof}

\begin{corollary} 
Soit $V_{\bullet}$ la filtration de Malgrange-Kashiwara de $G$ le long de $\tau =0$. Alors
$$F_{p}H_{\alpha}=(V_{\alpha}G\cap\tau^{p}G_{0}+V_{<\alpha}G)/(V_{<\alpha}G).$$
\end{corollary}
\qed\\

\noindent{Dans} tout ce qui suit, nous utiliserons indiff\'erement la notation $V_{\bullet}$ et ${\cal N}_{\bullet}$.
Nous venons de voir que les $F_{p}H_{\alpha}$ sont des $\cit [x]$-modules munis d'une connexion.  
Le lemme qui suit, dont la preuve est \'evidente, en donne des sections horizontales.
 Soit $\overline{\omega}_{p,\alpha}\in F_{p}H_{\alpha}$ : 
$\overline{\omega}_{p,\alpha}$ est l'image, par projection, d'un \'el\'ement 
$\omega_{\alpha +p}\tau^{p}$ avec 
$\omega_{\alpha +p}\in {\cal N}_{\alpha +p}\Omega^{n}(U)[x]$.

\begin{lemma} $\nabla_{\partial_{x_{i}}} \overline{\omega}_{p,\alpha}=0$ si et seulement si la classe de 
$\partial_{x_{i}}(\omega_{\alpha +p})$ appartient \`a $V_{<\alpha +p}G\cap G_{0}$. 
En particulier, $\overline{\omega}_{p,\alpha}$ est une section
horizontale de $\nabla$ si $\omega_{\alpha +p}$ ne d\'epend pas de $x$.
\end{lemma}

\begin{corollary}
Toute forme diff\'erentielle $\omega_{\alpha +p}\in {\cal N}_{\alpha +p}\Omega^{n}(U)$ d\'efinit, par projection, une section
horizontale $\overline{\omega}_{p,\alpha}$
de $(F_{p}H_{\alpha},\nabla )$ qui, pour la valeur nulle des param\`etres, coincide avec $\overline{\omega}_{p,\alpha}^{0}$, l'image de
 $\omega_{\alpha +p}$ dans $F_{p}H_{\alpha}^{0}$. Une telle section est unique. 
\end{corollary}

\noindent{Par} construction, toute base sur $\cit$ de $F_{p}H_{\alpha}^{0}$ fournit de cette mani\`ere
un syst\`eme de g\'en\'erateurs et  une base de $F_{p}H_{\alpha}$ sur $\cit [x]$.

\section{Le probl\`eme de Birkhoff en famille : le cas sous-diagramme}

Dans ce paragraphe, on garde les notations et les hypoth\`eses du d\'ebut du paragraphe 2.3. On a donc
$$F(u,x)=f(u)+\sum_{j=1}^{r}x_{j}g_{j}(u)$$
\noindent avec $\phi (g_{j})<1$ pour tout $j\in \{1,\cdots ,r\}$.

\subsection{$V$-solutions}
Le probl\`eme de Birkhoff pour $G_{0}$ est le suivant : trouver une base $\epsilon$ de $G_{0}$ sur $\cit [x,\theta ]$ dans laquelle la matrice de 
$\theta^{2}\nabla_{\partial_{\theta}}$ s'\'ecrit 
$$A_{0}(x)+A_{1}(x)\theta$$
et celle de $\nabla_{\partial_{x_{i}}}$ s'\'ecrit 
$$C^{(i)}_{-1}(x)\theta^{-1}+C^{(i)}_{0}(x)$$
pour tout $i\in \{1,\cdots ,r\}$, les matrices
$A_{0}(x)$, $A_{1}(x)$, $C^{(i)}_{-1}(x)$ et $C^{(i)}_{0}(x)$ \'etant \`a coefficients dans $\cit [x]$.

\begin{proposition} On suppose qu'il existe un $\cit [x,\tau ]$-sous-module libre $G_{\infty}$ de $G$, de rang maximum, 
v\'erifiant les conditions suivantes :\\

{\bf a.} pour tout $\alpha$, on a une d\'ecomposition
$${\cal N}_{\alpha}G\cap G_{0}={\cal N}_{\alpha}G\cap G_{0}\cap G_{\infty}\oplus 
\theta [{\cal N}_{\alpha -1}G\cap G_{0}],$$

{\bf b.} $\nabla_{\partial_{x_{i}}}G_{\infty}\subset G_{\infty}$ pour tout $i\in \{1,\cdots ,r\}$,\\

{\bf c.} $G_{\infty}$ est stable par $\tau\nabla_{\partial_{\tau}}$.\\

\noindent{Alors} le probl\`eme de Birkhoff pour $G_{0}$ a une solution, appel\'ee $V$-solution.
\end{proposition}

\noindent{Notons} que les conditions {\bf b.} et {\bf c.} peuvent \^etre unifi\'ees en demandant que la connexion de Gauss-Manin $\nabla$ soit \`a p\^oles
logarithmiques le long de $\tau =0$. Le point {\bf a.}, joint au point 3. de la proposition 2.3.1, montre en particulier que, pour tout $\alpha$, 
$\cal{N}_{\alpha}G\cap G_{0}\cap G_{\infty}$ est {\em libre} sur $\cit [x]$ (car projectif).\\

\begin{proof} Pour simplifier les notations, nous \'ecrirons $\partial_{\theta}$ 
pour $\nabla_{\partial_{\theta}}$ etc.. Choisissons 
$$J_{\alpha}\subset {\cal N}_{\alpha}G\cap G_{0}\cap G_{\infty},$$
\noindent isomorphe par projection \`a 
$$\frac{{\cal N}_{\alpha}G\cap G_{0}\cap G_{\infty}}{{\cal N}_{<\alpha}G\cap G_{0}\cap G_{\infty}},$$
\noindent ce dernier module est aussi libre sur $\cit [x]$ d'apr\`es la proposition 2.3.1.
 La d\'ecomposition {\bf a.} montre que toute base de $\oplus_{\alpha}J_{\alpha}$ sur
$\cit [x]$ fournit un syst\`eme de g\'en\'erateur de $G_{0}$ sur $\cit [x,\theta ]$. Dans toute base de ce type, la matrice de 
$\theta^{2}\partial_{\theta}$ a la forme attendue parce que $G_{\infty}$ est stable par $\tau\partial_{\tau}$. 
On a en fait
$$(\tau\partial_{\tau}+\alpha )G_{0}\cap G_{\infty}\cap V_{\alpha}G\subset G_{0}\cap G_{\infty}\cap V_{\alpha}G\oplus 
\tau (G_{0}\cap G_{\infty}\cap V_{\alpha +1}G).$$
Rappelons que 
$\nabla_{\partial_{x_{i}}}{\cal N}_{\alpha}G_{0}\subset \tau {\cal N}_{\alpha +1}G_{0}$. Ainsi, avec la d\'ecomposition {\bf a.}, on obtient
$$\nabla_{\partial_{x_{i}}}{\cal N}_{\alpha}G_{0}\subset \tau [{\cal N}_{\alpha +1}G\cap G_{0}\cap G_{\infty}]\oplus 
{\cal N}_{\alpha}G_{0}.$$
Si de plus $\nabla_{\partial_{x_{i}}}G_{\infty}\subset G_{\infty}$, on en d\'eduit 
$$\nabla_{\partial_{x_{i}}}({\cal N}_{\alpha}G\cap G_{0}\cap G_{\infty})\subset\tau ({\cal N}_{\alpha +1}G\cap 
G_{0}\cap G_{\infty})\oplus {\cal N}_{\alpha}G\cap G_{0}\cap G_{\infty}$$
et la proposition s'ensuit. 
\end{proof}

\noindent{Ainsi}, pour r\'esoudre le probl\`eme de Birkhoff pour $G_{0}$ il suffit de montrer l'existence de 
r\'eseaux $G_{\infty}$ v\'erifiant les conditions {\bf a.}, {\bf b.} et {\bf c.} de la proposition 3.1.1. 
L'id\'ee, maintenant \'eprouv\'ee, est de se ramener \`a un probl\`eme d'alg\`ebre lin\'eaire {\em via} une filtration 
oppos\'ee \`a le filtration $F_{\bullet}$.
 Rappelons que pour la valeur nulle des param\`etres, 
il existe, sur l'espace vectoriel $H^{o}_{\alpha}$, 
une filtration d\'ecroissante $U^{\bullet}$, 
stable par monodromie, c'est \`a dire $N(U^{p}H^{o}_{\alpha})\subset U^{p}H^{o}_{\alpha}$ o\`u
$N$ d\'esigne l'op\'erateur nilpotent induit par $\tau\partial_{\tau}+\alpha$, et oppos\'ee \`a la filtration de Hodge 
$F_{\bullet}H_{\alpha}^{o}$ : ceci r\'esulte du fait que $F_{\bullet}H_{\alpha}^{o}$ est la filtration de Hodge d'une structure de Hodge mixte d'apr\`es [Sab2].
On a donc une d\'ecomposition en somme directe
$$F_{p}H^{o}_{\alpha}=\oplus_{q\leq p}G^{o}_{\alpha ,q}$$
o\`u  $G^{0}_{\alpha ,q}=F_{q}H^{o}_{\alpha}\cap U^{q}H^{o}_{\alpha}$ (ce sont des espaces vectoriels de dimension finie) et on en 
tire dans ce cas pr\'ecis (sans param\`etres, donc) une solution au probl\`eme de Birkhoff (voir [DoSa, Appendix B]) : pr\'ecis\'ement, on en d\'eduit
 un $\cit [\tau ]$-r\'eseau $G_{\infty}^{o}$ du syst\`eme de Gauss-Manin $G^{o}$ 
 satisfaisant les conditions {\bf a.} et {\bf c.} de la proposition 3.1.1. 
Ceci \'etant, on obtient, par transport parall\`ele (corollaire 2.3.6), 
une d\'ecomposition de $F_{p}H_{\alpha}$ en somme directe
$$F_{p}H_{\alpha}=\oplus_{q\leq p}G_{\alpha ,q}$$
o\`u maintenant les $G_{\alpha ,q}$ sont des $\cit [x]$-modules libres. On {\em d\'efinit} alors 
$$U^{p}H_{\alpha}=\oplus_{q\geq p}G_{\alpha ,q}.$$
\noindent{On} a par construction (horizontalit\'e) 
$$\nabla_{\partial_{x_{i}}}(U^{p}H_{\alpha})\subset U^{p}H_{\alpha}.$$
 De plus 
$$N(U^{p}H_{\alpha})\subset U^{p}H_{\alpha}$$ 
parce qu'une telle relation est vraie pour la valeur nulle des param\`etres et que $N$ commute avec $\nabla_{\partial_{x_{i}}}$ et donc transforme sections horizontales en sections horizontales.

\begin{proposition} Le probl\`eme de Birkhoff pour $G_{0}$ a des solutions : il existe des $\cit [x,\tau ]$-sous-modules $G_{\infty}$
de $G$ 
v\'erifiant les conditions de la proposition 3.1.1.
\end{proposition}
\begin{proof}
Comme dans le cas absolu,
une filtration d\'ecroissante de $H_{\alpha}$, oppos\'ee \`a la filtration $F_{\bullet}$,
 stable par $N$ et  $\nabla_{\partial_{x_{i}}}$ pour tout $i$, fournit un
 r\'eseau $G_{\infty}$ v\'erifiant les conditions {\bf b.} et {\bf c.}
de la proposition 3.1.1 (ceci est une cons\'equence de [Sab3, premi\`ere partie de la preuve du th\'eor\`eme III.1.1] : voir [D, 6.2.2] pour une explication d\'etaill\'ee; 
la correspondance {\em bijective} entre deux tels objets reste \`a discuter, m\^eme si elle n'est pas n\'ecessaire ici). On a aussi 
la propri\'et\'e {\bf a.},
parce qu'elle est vraie pour toutes les valeurs des param\`etres et que les $\cit [x]$-modules apparaissant dans la d\'ecomposition
 sont tous libres de type fini. 
\end{proof}

\noindent{Remarquons} que, par construction, la restriction de $G_{\infty}$ \`a la valeur nulle des param\`etres coincide avec le r\'eseau $G_{\infty}^{o}$ 
consid\'er\'e dans le cas absolu. De plus, si l'on part de la filtration oppos\'ee canonique $U^{\bullet}$ de M. Saito [Sai2, lemma 2.8] 
sur $H^{o}_{\alpha}$, pour toute
valeur $x^{o}$ du param\`etre la restriction $G_{\infty}^{x^{o}}$ de $G_{\infty}$ \`a $x^{o}$ est la $V$-solution canonique du probl\`eme de Birkhoff 
pour le polyn\^ome $F(u,x^{o})$ consid\'er\'ee dans [DoSa, 3.c et appendix B].
Gardons les notations du d\'ebut de ce paragraphe. 
Les relations d'int\'egrabilit\'e fournissent alors le

\begin{corollary} 
Il existe une base $\epsilon$ de $G_{0}$ dans laquelle la matrice de $\theta^{2}\partial_{\theta}$ s'\'ecrit 
$$B_{0}(x)+B_{\infty}\theta$$
avec $B_{\infty}$ {\em constante} et celle de $\nabla_{\partial_{x_{i}}}$ s'\'ecrit
$$C^{(i)}(x)\theta^{-1}$$
$(i\in \{1,\cdots ,r\})$. Les matrices $B_{0}(x)$ et $C^{(i)}(x)$ sont \`a coefficients dans $\cit [x]$.
\end{corollary}

\begin{proof} Soit $\nu$ une base de $G_{0}$ construite dans la preuve de la proposition 3.1.1, $A_{0}(x)$, $A_{1}(x)$, $C^{(i)}_{0}(x)$ et $C^{(i)}_{-1}(x)$
les matrices correspondantes (voir d\'ebut du paragraphe 3.1). Les matrices $C^{(i)}_{0}(x)$ v\'erifient, parce que la connexion est int\'egrable,
$$\frac{\partial }{\partial x_{i}}C^{(j)}_{0}(x)-\frac{\partial }{\partial x_{j}}C^{(i)}_{0}(x)=[C^{(j)}_{0}(x),C^{(i)}_{0}(x)]$$
pour tous $i,j$. Il existe donc une unique matrice $P$ \`a coefficients {\em a priori} holomorphes telle que $P_{|x=0}=Id$ et 
$\frac{\partial }{\partial x_{i}}P=-C^{(i)}_{0}(x)P$. 
Mais, comme on a, par hypoth\`ese sur les $g_{i}$,
$$g_{i}\cal{N}_{\alpha}K[x]\subset\cal{N}_{<\alpha +1}K[x],$$
\noindent on peut supposer 
(il suffit d'ordonner les bases selon le 
poids croissant) que les matrices $C^{(i)}_{0}(x)$ sont nilpotentes, triangulaires. Nous sommes donc r\'eduits \`a r\'esoudre un syst\`eme, en l'inconnue $y$,
du type
$$\frac{\partial }{\partial x_{k}}y=b_{k}(x),$$
\noindent avec $b_{k}(x)\in\cit [x]$. $y$ est polynomiale et la matrice $P$ est ainsi \`a coefficients polynomiaux. Pour terminer la preuve du 
corollaire, il suffit de poser $\epsilon =P\nu$ et de remarquer que $P^{-1}$ est aussi \`a coefficients polynomiaux pour les m\^emes raisons que ci-dessus. 
\end{proof} 

\noindent Par la suite, nous ne consid\'ererons que des bases donn\'ees par les corollaire 3.1.3. Nous les appellerons {\em bases plates} (voir le paragraphe 4 pour une justification de cette terminologie).
Signalons que les matrices $B_{0}$, $B_{\infty}$ et $C^{(i)}$ fournies par le corollaire 3.1.3 v\'erifient les relations d'int\'egrabilit\'e\\

{\bf I.1}   $\partial_{x_{j}}C^{(i)}=\partial_{x_{i}}C^{(j)}$, pour $i,j=1,\cdots ,r$\\

{\bf I.2}   $[C^{(i)},C^{(j)}]=0$, pour $i,j=1,\cdots ,r$\\
 
{\bf I.3}   $[B_{0},C^{(i)}]=0$, pour $i=1,\cdots ,r$\\

{\bf I.4}   $\partial_{x_{i}}B_{0}+C^{(i)}=[B_{\infty},C^{(i)}]$, pour $i=1,\cdots ,r$.\\

\noindent{Par} unicit\'e, la base $\epsilon$ donn\'ee par le corollaire 3.1.3 coincide avec celle fournie par le
th\'eor\`eme g\'en\'eral de B. Malgrange [Mal, th\'eor\`eme 2.2]. 
En fait, on a ici un peu plus : la restriction de la solution $\epsilon$ \`a une valeur $x^{o}$ quelconque du param\`etre est une solution 
du probl\`eme de Birkhoff donn\'ee par [DoSa] (voir aussi la preuve de la proposition 3.1.1) 
pour le polyn\^ome $F(u,x^{o})$ : ceci explique l'assertion de l'introduction.

\subsection{$V^{+}$-solutions}

Si l'on part d'une filtration $U^{\bullet}$ oppos\'ee \`a la filtration $F_{\bullet}$ sur chaque $H_{\alpha}^{o}$ et v\'erifiant de plus 
$N(U^{p}H_{\alpha}^{o})\subset U^{p+1}H_{\alpha}^{o}$ (c'est le cas en particulier si $U^{\bullet}$ est la filtration oppos\'ee canonique de M. Saito)
 on obtient comme ci-dessus un r\'eseau $G_{\infty}$ de $G$ v\'erifiant les conditions {\bf a.}, 
{\bf b.} et {\bf c.} de la proposition 3.1.1 avec de plus
$$(\tau\partial_{\tau}+\alpha )G_{0}\cap G_{\infty}\cap V_{\alpha}G\subset G_{0}\cap G_{\infty}\cap V_{<\alpha}G\oplus 
\tau (G_{0}\cap G_{\infty}\cap V_{\alpha +1}G).$$
$G_{\infty}$ est une $V^{+}$-solution au sens de [DoSa, appendix B].
Les conclusions du corollaire 3.1.3 restent vraies, bien entendu : la matrice $B_{\infty}$ est maintenant {\em diagonalisable}. Par construction, l'ensemble de ses valeurs propres coincide avec le spectre de $(G^{o}, G_{0}^{o})$ ou encore avec le spectre \`a l'infini du polyn\^ome de Laurent $f$.\\

\subsection{Dualit\'e}

Le but de ce paragraphe est de montrer comment, dans notre situation, la dualit\'e du cas absolu se prolonge alg\'ebriquement en famille. 
Soit $\epsilon$ une base de $G_{0}$ dans laquelle la matrice de la connexion s'\'ecrit
$$(\frac{B_{0}(x)}{\theta}+B_{\infty})\frac{d\theta}{\theta}+\sum_{i=1}^{r}\frac{C^{(i)}(x)}{\theta}dx_{i}.$$
\noindent Une telle base existe d'apr\`es le corollaire 3.1.3 et on peut, d'apr\`es 3.2, supposer que la matrice $B_{\infty}$ est diagonale, 
$$B_{\infty}=diag\ (\alpha_{1},\cdots ,\alpha_{\mu}),$$ 
avec $\alpha_{1}\leq\alpha_{2}\leq\cdots \leq\alpha_{\mu}$.

\subsubsection{Dualit\'e absolue}
Pour chaque valeur $x$ fix\'ee des param\`etres, le syst\`eme de Gauss-Manin $G^{x}$ (restriction de $G$ \`a $x$) 
est muni d'une dualit\'e $S^{x}$, forme sesquilin\'eaire non d\'eg\'en\'er\'ee, compatible aux connexions, avatar de la dualit\'e de Poincar\'e
(voir par exemple [Sab3, chapitre III, 1.b] pour une d\'efinition g\'en\'erale de ces notions et [DoSa2, p. 9] pour une transposition 
\`a la situation g\'eom\'etrique consid\'er\'ee ici). Dans la situation du corollaire 3.1.3, et d'apr\`es [Sai2, lemma 2.8], on peut supposer que la base 
$\epsilon^{o}=(\epsilon_{1}^{o},\cdots ,\epsilon_{\mu}^{o})$, restriction de la base $\epsilon$ \`a la valeur nulle des param\`etres, est une $S^{o}$-solution orthonorm\'ee, c'est \`a dire que
$$S^{o}(\epsilon_{i}^{o},\epsilon_{\mu +1-j}^{o})=\delta_{ij}\tau^{-n}.$$
Pour tout $x\in\cit^{n}$, et pour la m\^eme raison d'apr\`es la remarque qui suit la proposition 3.1.2, $\epsilon^{x}$ est une 
$S^{x}$-solution, {\em i.e}
$$S^{x}(\epsilon_{i}^{x},\epsilon_{j}^{x})=S^{x}_{-n}(\epsilon_{i}^{x},\epsilon_{j}^{x})\tau^{-n}$$
avec $S^{x}_{-n}(\epsilon_{i}^{x},\epsilon_{j}^{x})\in\cit$. Ce coefficient ne d\'epend que de la classe de 
$\epsilon_{i}^{x}$ et $\epsilon_{j}^{x}$ dans $G_{0}^{x}/\theta G_{0}^{x}$ : c'est le classique r\'esidu de Grothendieck (se r\'ef\'erer par exemple
\`a [He, section 10.6] et aux r\'ef\'erences qui y sont mentionn\'ees). On a en fait un peu plus

\begin{lemma} Pour tout $x\in\cit^{n}$, on a $S^{x}_{-n}(\epsilon_{i}^{x},\epsilon_{j}^{x})=S^{o}_{-n}(\epsilon_{i}^{o},\epsilon_{j}^{o})$. En particulier,
$$S^{x}(\epsilon_{i}^{x},\epsilon_{j}^{x})=S^{o}(\epsilon_{i}^{o},\epsilon_{j}^{o}).$$
\end{lemma}
\begin{proof} L'outil essentiel est le comportement de la dualit\'e vis-\`a-vis de la $V$-filtration. Soit $\alpha_{i}$ l'ordre
de $\epsilon_{i}^{x}$ pour la $V$-filtration et 
$\epsilon_{i}^{x}\in\eta_{i}^{x}+V_{<\alpha_{i}}.$ 
La construction de la section 3.1 montre que les parties principales $\eta_{i}^{x}$ sont constantes ({\em i.e} ne d\'ependent pas de $x$) : 
on les notera donc $\eta_{i}$. 
Remarquons que $S^{x}_{-n}(\epsilon_{i}^{x},\epsilon_{j}^{x})=0$ si 
$\alpha_{i}+\alpha_{j}\neq n$. Soient donc $\epsilon_{i}^{x}$ et $\epsilon_{j}^{x}$ tels que
$\alpha_{i}+\alpha_{j}= n$.
On a alors $S^{x}_{-n}(\epsilon_{i}^{x},\epsilon_{j}^{x})=
S^{x}_{-n}(\eta_{i},\eta_{j})$ pour des raisons d'orthogonalit\'e (voir par exemple [He, Theorem 10.28]). 
Comme $\eta_{i}$ et $\eta_{j}$ sont homog\`enes pour la $V$-filtration, cette derni\`ere quantit\'e se calcule 
en consid\'erant la restriction du r\'esidu \`a 
$$gr_{V}^{\alpha }G_{0}^{x}/\theta G_{0}^{x}\times gr_{V}^{n-\alpha}G_{0}^{x}/\theta G_{0}^{x}.$$
Mais, toujours gr\^ace \`a l'hypoth\`ese faite sur les mon\^omes de la d\'eformation, on a
$$gr_{V}^{\beta }G_{0}^{x}/\theta G_{0}^{x}=gr_{V}^{\beta }G_{0}^{o}/\theta G_{0}^{o}$$
pour tout $\beta$ et ainsi 
$S^{x}_{-n}(\eta_{i},\eta_{j})=S^{o}_{-n}(\eta_{i},\eta_{j})$. Le lemme est montr\'e.
\end{proof}

\subsubsection{Sym\'etrie} 
On d\'efinit la transformation
$$T:\cal{M}_{\mu\times\mu}(\cit )\rightarrow \cal{M}_{\mu\times\mu}(\cit )$$
en posant $TA=(a_{\mu +1-j \ \mu +1-i})$ si $A=(a_{ij})$. Nous dirons qu'une matrice $A$ est $T$-sym\'etrique si elle v\'erifie $TA=A$. 

\begin{proposition} Dans la situation du corollaire 3.1.3, on a :\\ 
1. la matrice $B_{0}(x)$ est $T$-sym\'etrique pour tout $x\in\cit^{n}$ et $B_{\infty}+TB_{\infty}=nId$,\\
2. les matrices $C^{(i)}(x)$, $i\in \{1,\cdots ,r\}$, sont uniquement d\'etermin\'ees par la matrice $B_{0}(x)$ via la formule
$$\partial_{x_{i}}B_{0}(x)+C^{(i)}(x)=[B_{\infty} , C^{(i)}(x)],$$
3. les matrices $C^{(i)}(x)$, $i\in \{1,\cdots ,r\}$, sont $T$-sym\'etriques pour tout $x\in\cit^{n}$.
\end{proposition}
\begin{proof} 1. $B_{0}(x)$ est auto-adjointe pour $S^{x}$ parce que $S^{x}$ est compatible aux connexions. La propri\'et\'e de sym\'etrie attendue 
r\'esulte alors du lemme 3.3.1. M\^eme chose pour $B_{\infty}$.\\
2. Ecrivons $C^{(i)}(x)=(c^{i}_{k\ell})$, $B_{0}(x)=(b^{0}_{k\ell})$ et $B_{\infty}=diag(\alpha_{1},\cdots ,\alpha_{\mu})$. 
La formule \'evoqu\'ee dans l'\'enonc\'e est la relation d'int\'egrabilit\'e {\bf I.4}
de la connexion et s'\'ecrit
$$(1-\alpha_{k}+\alpha_{\ell})c^{i}_{k\ell}=\partial_{x_{i}}b^{0}_{k\ell}.$$
Si $\alpha_{k}=1+\alpha_{\ell}$ alors $c^{i}_{k\ell}=0$ parce que, avec nos hypoth\`eses (les $g_{i}$ 
sont 'strictement sous-diagramme'),
 $$g_{i}{\cal N}_{\alpha_{\ell}}K[x]\subset {\cal N}_{<\alpha_{\ell}+1}K[x].$$
\noindent Sinon, on a
$$c^{i}_{k\ell}=\frac{1}{1-\alpha_{k}+\alpha_{\ell}}\partial_{x_{i}}b^{0}_{k\ell}.$$\\
3. Par d\'efinition, $Tc^{i}_{k\ell}=c^{i}_{\mu +1-\ell\ \mu +1-k}$ et cette derni\`ere quantit\'e vaut 
$$\frac{1}{1-\alpha_{\mu +1-\ell}+\alpha_{\mu +1-k}}\partial_{x_{i}}b^{0}_{\mu +1-\ell\ \mu +1-k}$$ 
d'apr\`es 2. Il suffit maintenant d'utiliser le point 1.
 \end{proof}

\subsubsection{Dualit\'e en famille}
Toujours dans la situation du corollaire 3.1.3, on pose 
$$S(\epsilon_{i},\epsilon_{j})=S^{o}(\epsilon_{i}^{o},\epsilon_{j}^{o})$$
formule que l'on \'etend \`a $G$ par $\cit [x][\tau ,\tau^{-1}]$-(sesqui)lin\'earit\'e.
Rappelons que $S$ est dite de poids $w$ relativement au r\'eseau $G_{0}$ si
$$S(G_{0},G_{0})\subset \theta^{w}\cit [x,\theta ].$$

\begin{proposition} $S$ est une forme sesquilin\'eaire non d\'eg\'en\'er\'ee, compatible aux connexions, de poids $n$ relativement au r\'eseau $G_{0}$.
\end{proposition}
\begin{proof} Il suffit de remarquer que les matrices $B_{0}(x)$ et $C^{(i)}(x)$, $i\in \{1,\cdots ,r\}$, sont $T$-sym\'etriques. 
C'est pr\'ecis\'ement ce que donne la proposition 3.3.2.
\end{proof}  

\begin{remark} L'existence de cette dualit\'e $S$ devrait pouvoir aussi se montrer en invoquant les arguments g\'en\'eraux
(dualit\'e des ${\cal D}$-modules) donn\'es dans [DoSa, section 2.b] et [Sab1, section 11] : la compatibilit\'e par rapport aux connexions montre directement que $S$ est constante lorsqu'elle est \'evalu\'ee sur la base $\epsilon$.
\end{remark}

\section{Application \`a la construction de vari\'et\'es de Frobenius}
\subsection{Vari\'et\'es de Frobenius et polyn\^omes de Laurent}
Soit $f:U\rightarrow\cit$ un polyn\^ome de Laurent commode et non d\'eg\'en\'er\'e, $\mu$ son nombre de Milnor global et $\tilde{M}$ un voisinage 
ouvert de l'origine dans $\cit^{\mu}$. On d\'efinit
$$\tilde{F}:(\cit^{*})^{n}\times \tilde{M}\rightarrow \cit$$
par 
$$\tilde{F}(u,x)=f(u)+\sum_{i=1}^{\mu}x_{i}g_{i}(u)$$
 o\`u les $g_{i}$ sont des polyn\^omes de Laurent se projetant en une base de $\cit [u,u^{-1}]/(\partial_{u_{1}} f,\cdots ,\partial_{u_{n}}f)$. 
La th\'eorie g\'en\'erale, expos\'ee dans [DoSa], nous apprend que l'on peut munir $\tilde{M}$ d'une structure de Frobenius. Le but de ce paragraphe est de donner une preuve plus simple de ce fait 
en utilisant [HeMa, theorem 4.5] et les r\'esultats 
des paragraphes 2 et 3 de ces notes. On \'evite en 
particulier le d\'elicat passage analytique et le recours aux proc\'ed\'es d'alg\'ebrisation transcendants utilis\'es dans les paragraphes 1 et 2 de [DoSa].
L'id\'ee est que l'on peut, sous certaines conditions, d\'ecrire la structure de Frobenius ainsi construite sur $\tilde{M}$ \`a partir de donn\'ees sur un sous-espace $M$ de 
$\tilde{M}$ : $M$ sera un voisinage ouvert de l'origine dans $\cit^{r}$, espace des param\`etres d'une d\'eformation $F$ de $f$ d\'efinie par
$$F(u,x)=f(u)+\sum_{j=1}^{r}x_{j}g_{j}(u).$$
L'espoir est qu'il existe une d\'eformation de ce type, {\em avec $\phi (g_{j})<1$ pour tout $j\in \{1,\cdots ,r\}$}, qui v\'erifie les hypoth\`eses du th\'eor\`eme de Manin et Hertling.\\

\noindent Avant de d\'etailler ce point, r\'esumons les r\'esultats des paragraphes 2 et 3 obtenus pour de telles d\'eformations : au paragraphe 2, nous avons montr\'e que le r\'eseau de Brieskorn associ\'e est libre, de rang $\mu$, sur $\cit [x,\theta ]$.
Au paragraphe 3 nous avons ensuite r\'esolu 
le probl\`eme de Birkhoff pour $G_{0}$ : ceci signifie que l'on peut \'etendre $G_{0}$ en un fibr\'e trivial ${\cal G}$ 
sur $\ppit^{1}\times M$, muni d'une connexion m\'eromorphe int\'egrable $\nabla$ \`a p\^oles logarithmiques le long de $\{\infty \}\times M$ et d'ordre $1$ le long de $\{0\}\times M$. On obtient ainsi, en consid\'erant de plus la dualit\'e d\'efinie au paragraphe 3.3, une $trTLEP(n)$-structure au sens de [HeMa, definition 4.1]. Cette $trTLEP(n)$-structure permet \`a son tour (voir [HeMa, theorem 4.2] et les r\'ef\'erences qui y sont mentionn\'ees) de d\'efinir une structure de type Frobenius 
$$(M,E, \bigtriangledown , \Phi , R_{0}, R_{\infty}, g)$$ au sens de [HeMa, definition 4.1]. 
On a ici $E={\cal G}_{|\{0\}\times M}$.
Avec les notations de 3.3, $g:=S_{-n}$ et la connexion $\bigtriangledown$ d\'efinie dans cet uplet est la connexion de Levi-Civita de $g$. Sa matrice est nulle dans la base $[\epsilon ]$ de $E$ induite par la solution du probl\`eme de Birkhoff $\epsilon$. $R_{0}$ ({\em resp.} $R_{\infty}$) est l'endomorphisme de $E$ de matrice $B_{0}$ ({\em resp.} $B_{\infty}$) dans la base $[\epsilon ]$ et $\Phi :=\sum_{i}\Phi^{(i)}dx_{i}$ est l'application de $E$ dans $\Omega^{1}(M)\otimes E$ d\'efinie par $\sum_{i}C^{(i)}dx_{i}$ dans la base $[\epsilon ]$. $\Phi^{(i)}$ est ainsi la multiplication par $-g_{i}$.
\noindent Notons 
$$E_{0}:=\frac{\Omega^{n}(U)}{df\wedge\Omega^{n-1}(U)}.$$ 
On d\'esignera par les m\^emes lettres la restriction des objets $\Phi$, $R_{0}$, $R_{\infty}$... \`a $E_{0}$.
Le th\'eor\`eme 4.5 de [HeMa] affirme que, moyennant les conditions rappel\'ees ci-dessous, 
la $trTLEP(n)$-structure obtenue ci-dessus sur $M$ se d\'eploie de mani\`ere unique en une $trTLEP(n)$-structure sur $\tilde{M}$. Une forme primitive et homog\`ene permet ensuite, {\em via} l'application de p\'eriodes associ\'ee, de montrer que $\tilde{M}$ est une vari\'et\'e de Frobenius.\\

\noindent Les conditions figurant dans le th\'eor\`eme 4.5 de [HeMa] sont les suivantes:\\

\noindent $\bullet$ {\em La condition (EC).---} Elle demande l'existence d'un \'el\'ement $\zeta\in E_{0}$, vecteur propre de $R_{\infty}$ et v\'erifiant les conditions \'enum\'er\'ees ci-dessous. Un choix naturel, et nous verrons pourquoi, est de prendre pour $\zeta$ la classe de la forme volume $du/u:=du_{1}/u_{1}\wedge\cdots\wedge du_{n}/u_{n}$. 
Comme $f$ est un polyn\^ome de Laurent commode et non d\'eg\'en\'er\'e, la forme volume $du/u:=du_{1}/u_{1}\wedge\cdots\wedge du_{n}/u_{n}$ induit bien un vecteur propre de $R_{\infty}$
(voir [DoSa, 4.d] et [D, proposition 7.0.2]).\\

\noindent {\em Dans tout ce qui suit, $\zeta$ d\'esignera la classe de $du/u $ dans $E_{0}$}.\\

\noindent $\bullet$ {\em La condition (IC).---}  Elle se lit :
\begin{center} {\em l'application $\Phi_{\bullet}:T_{0}M\rightarrow E_{0}$, 
$X\mapsto \Phi_{X}\zeta$ est injective.}
\end{center}
Ici, $\Phi_{X}$ d\'esigne l'endomorphisme obtenu par contraction de $\Phi$ avec $X$. On a ainsi, avec les notations pr\'ec\'edentes, $\Phi_{\partial_{x_{i}}}=\Phi^{(i)}$. Rappelons que, par d\'efinition, $\Phi^{(i)}(\zeta)$ est \'egal \`a la classe de $-g_{i}\frac{du}{u}$ dans $E_{0}$.
\begin{lemma} On suppose que les polyn\^omes de Laurent $g_{j}$ sont lin\'eairement ind\'ependants dans $\cit [u,u^{-1}]$. Si de plus $\phi (g_{j})<1$ pour tout $j\in \{1,\cdots ,r\}$, les classes de $\Phi^{(1)}(\zeta ),\cdots ,\Phi^{(r)}(\zeta )$ dans $E_{0}$ sont lin\'eairement ind\'ependantes. 
\end{lemma}
\begin{proof} Il suffit de v\'erifier que les classes de $g_{1},\cdots ,g_{r}$ dans
$\cit [u,u^{-1}]/(\partial_{u_{1}}f,\cdots ,\partial_{u_{n}}f)$ sont ind\'ependantes. Mais ceci r\'esulte des conditions 
$\phi (g_{j})<1$ : en effet, supposons qu'il existe des nombres complexes $\alpha_{1},\cdots ,\alpha_{r}$ tels que
$$\sum_{j=1}^{r}\alpha_{j}g_{j}=\sum_{i=1}^{n}b_{i}u_{i}\frac{\partial f}{\partial u_{i}}.$$
On peut choisir, d'apr\`es la formule de division \'enonc\'ee dans la preuve du lemme 2.1.1, les $b_{i}$ tels que $\phi (b_{i})\leq \phi (\sum_{j=1}^{r}\alpha_{j}g_{j})-1$ donc $\phi (b_{i})<0$ parce que $\phi (g_{j})<1$ pour tout $j$. Ainsi, $b_{i}=0$ pour tout $i$. Si les polyn\^omes $g_{j}$ sont ind\'ependants dans $\cit [u,u^{-1}]$, ceci exige que $\alpha_{i}=0$ pour tout $i$.
\end{proof}

\noindent Ainsi, la condition (IC) est v\'erifi\'ee lorsque les hypoth\`eses du lemme 4.1.1 sont satisfaites.\\

\noindent $\bullet$ {\em La condition (GC).---}  C'est cette condition qui nous permettra de fixer d\'efinitivement la d\'eformation $F$. Elle se lit : 
\begin{center}
{\em $\zeta$ et ses images par l'it\'eration des applications $R_{0}$ et $\Phi^{(i)}$, 
$i\in \{1,\cdots ,r\}$,  engendrent $E_{0}$}. 
\end{center}
Bien s\^ur, cette condition est
en particulier v\'erifi\'ee si 
$\zeta$ et ses images par l'it\'eration des applications $\Phi^{(i)}$ engendrent $E_{0}$. Cette remarque nous permet d'\'enoncer (rappelons que $\cit [g_{1},\cdots ,g_{r}]$ d\'esigne la sous alg\`ebre des endomorphismes de $E_{0}$ engendr\'ee par la multiplication par les $g_{i}$), \`a l'aide de [HeMa, theorem 4.5] et de sa preuve, le

\begin{theorem} Soit $f$ un polyn\^ome de Laurent non d\'eg\'en\'er\'e et commode, $g_{1},\cdots ,g_{r}$ des \'el\'ements ind\'ependants de $\cit [u,u^{-1}]$ tels que $\phi (g_{j})<1$ pour tout $j\in \{1,\cdots ,r\}$. On suppose de plus que $\cit [g_{1},\cdots ,g_{r}]\zeta =E_{0}$. Alors,\\
$1.$ la trTLEP(n)-structure fournie, selon le proc\'ed\'e du paragraphe 3, par le r\'eseau de Brieskorn de la d\'eformation 
$$F(u,x)=f(u)+\sum_{j=1}^{r}x_{j}g_{j}(u)$$
se d\'eploie de mani\`ere universelle,\\
$2.$ $\tilde{M}$ est une vari\'et\'e de Frobenius.
\end{theorem}

\begin{example} On suppose que $\phi (u_{i})<1$ et $\phi (1/u_{i})<1$ pour tout $i\in\{1,\cdots ,n\}$. Alors la construction pr\'ec\'edente, avec la d\'eformation $$F(u,x)=f(u)+\sum_{i=1}^{n}x_{i}u_{i}+\sum_{i=1}^{n}x_{n+i}\frac{1}{u_{i}},$$
permet toujours de munir $\tilde{M}$ d'une structure de Frobenius.
\end{example}

\begin{remark} $1.$ La structure de Frobenius obtenue sur $\tilde{M}$ d\'epend donc de deux choix : celui de la forme homog\`ene $\zeta$ et celui de la d\'eformation $F$.\\ 
$2.$ La condition (GC) est aussi v\'erifi\'ee si $\zeta$ et ses images par l'it\'eration de l'applications $R_{0}$ engendrent $E_{0}$. On peut alors prendre $M=\{point \}$ et $F=f$ 
(la condition (IC) est dans ce cas vide) : c'est la situation originellement trait\'e par B. Malgrange. 
Ceci se produit par exemple si toutes les valeurs 
critiques de $f$ sont distinctes : la condition (GC) est v\'erifi\'ee parce que les valeurs propres de $B_{0}$ (c'est la matrice de la multiplication par $f$ sur $E_{0}$ dans la base $[\epsilon ]$) sont toutes distinctes. La classe d'exemples \'etudi\'ee dans [DoSa2] rentre dans ce cadre. 
\end{remark}

\subsection{Vari\'et\'es de Frobenius et polyn\^omes}
\subsubsection{}
La construction pr\'ec\'edente peut aussi se faire, moyennant de l\'eg\`eres modifications,
lorsque $U=\cit^{n}$ et 
$$f:U\rightarrow\cit$$
\noindent est une application polynomiale, commode et non d\'eg\'en\'er\'ee. 
$f$ est alors \`a points critiques isol\'es et son nombre de Milnor (global) est d\'efini par la formule
$$\mu =\dim_{\cit}\Omega^{n}(U)/df\wedge \Omega^{n-1}(U),$$
$\Omega^{p}(U)$ d\'esignant l'espace des $p$-formes r\'eguli\`eres sur $U$. Si l'on garde \`a l'esprit la condition (GC), la d\'eformation de $f$ naturelle \`a consid\'erer est
$$F(u,x)=f(u)+\sum_{i=1}^{n}x_{i}u_{i}.$$
On d\'efinit, de la m\^eme mani\`ere qu'au paragraphe 2, le module de Brieskorn $G_{0}$ et le syst\`eme de Gauss-Manin $G$ associ\'es \`a $F$ et, comme en 2.1.1, une filtration de Newton $\phi$ sur $K=\cit [u_{1},\cdots ,u_{n}]$.  

\begin{proposition}  On suppose que $\phi (u_{j})<1$ pour tout $j\in \{1,\cdots ,n\}$. Alors
tous les r\'esultats des paragraphes 2 et 3 restent valides pour la d\'eformation
$$F(u,x)=f(u)+\sum_{j=1}^{n}x_{j}u_{j}$$
du polyn\^ome $f$. En particulier $G_{0}$ est libre de rang $\mu$ sur $\cit [x,\theta ]$ et le probl\`eme de Birkhoff en famille admet une solution alg\'ebrique.
\end{proposition}

\noindent Avant de montrer la proposition 4.2.1, remarquons que les r\'esultats du paragraphe 2 
s'adaptent comme suit : si $K=\cit [u_{1},\cdots ,u_{n}]$, on d\'efinit 
$${\cal N}_{\alpha}K[x]:=\{v\in K[x], \phi ^{*}(v)\leq\alpha \}$$
\noindent o\`u $\phi ^{*}(v)=\phi (u_{1}\cdots u_{n}v)$. On en d\'eduit une filtration ${\cal N}_{\bullet}$ sur $\Omega^{n}(U)[x]$ en posant
$${\cal N}_{\alpha}\Omega^{n}(U)[x]:=\{gdu\in \Omega^{n}(U)[x], \phi^{*} (g)\leq\alpha \}$$
\noindent puis sur $\Omega^{n}(U)[x,\theta ]$, $G_{0}$ et $G$ de la m\^eme mani\`ere qu'en 2.1.
On d\'efinit de m\^eme les modules $E_{\beta}[x]$ pour $\beta\in\qit$.
Soit
$$F(u,x)=f(u)+\sum_{i=1}^{r}x_{i}g_{i},$$
\noindent les $g_{j}$ v\'erifiant  
$$\phi (\frac{\partial g_{j}}{\partial u_{i}})<1-\phi (u_{i})$$
\noindent pour tout $i\in \{1,\cdots, n\}$ et tout $j\in \{1,\cdots ,r\}$. Cette condition doit \^etre compar\'ee avec celle donn\'ee au lemme 2.1.1 dans le cas des polyn\^omes de Laurent. Il faut aussi noter  que la relation entre les conditions 
$\phi (\frac{\partial g_{j}}{\partial u_{i}})<1-\phi (u_{i})$
et $\phi (g_{j})<1$ n'est en g\'en\'eral pas tr\`es  claire et ceci expliquera pourquoi la proposition 4.2.1 n'est {\em a priori} valable que pour les d\'eformations qui y sont mentionn\'ees. 

\begin{lemma} Sous les hypoth\`eses pr\'ec\'edentes, il existe, 
pour tout $h\in{\cal N}_{\alpha}K[x]$, $v_{\alpha}\in \oplus_{\beta\leq\alpha}E_{\beta}[x]$ et $a_{1},\cdots a_{n}\in K[x]$ tels que
$$h=v_{\alpha}+\sum_{i=1}^{n}a_{i}\frac{\partial F}{\partial u_{i}}$$
les $a_{i}$ v\'erifiant de plus $\phi^{*} (a_{i})\leq \alpha -1+\phi (u_{i})$ et
$\phi^{*} (\partial a_{i}/\partial u_{i})\leq \alpha -1$. 
\end{lemma} 
\begin{proof} Analogue \`a celle du lemme 2.1.1. en utilisant le th\'eor\`eme de division 
dans le cas sans param\`etres qui s'\'ecrit maintenant
$$h=v_{\alpha}+\sum_{i=1}^{n}b_{i}\frac{\partial f}{\partial u_{i}}$$
\noindent avec $\phi^{*} (b_{i})\leq \alpha -1+\phi (u_{i})$ et $\phi^{*} (\partial b_{i}/\partial u_{i})\leq \alpha -1$.
\end{proof} 

\noindent {\em Preuve de la preuve de la proposition 4.2.1} : on prend ici 
$n=r$ et $g_{i}=u_{i}$ pour tout $i$. La condition $\phi (\frac{\partial u_{j}}{\partial u_{i}})<1-\phi (u_{i})$ s'\'ecrit alors, et c'est un petit miracle,
$$\phi (u_{i})<1.$$
\noindent L'in\'egalit\'e  $\phi (\frac{\partial u_{j}}{\partial u_{i}})<1-\phi (u_{i})$ permet, {\em via} le lemme 4.2.2, d'obtenir les r\'esultats du paragraphe 2. L'in\'egalit\'e $\phi (u_{i})<1$ permet d'obtenir, avec les m\^emes preuves, les r\'esultats du paragraphe 3.\qed \\

\noindent Concernant le paragraphe 4, et donc la construction de vari\'et\'es de Frobenius, deux 
nouvelles difficult\'es se pr\'esentent, la premi\`ere \'etant s\'erieuse :\\

\noindent $\bullet$ comme plus haut, un choix naturel pour la
forme $\zeta$ v\'erifiant les conditions du th\'eor\`eme de Manin et Hertling 
est la classe de la forme volume $du_{1}\wedge\cdots \wedge du_{n}$. Dans notre situation, son homog\'en\'eit\'e n'est pas acquise en toute g\'en\'eralit\'e. Elle l'est cependant si 
$$du_{1}\wedge\cdots \wedge du_{n}\in {\cal N}_{\alpha_{1}}\Omega^{n}(U),$$
\noindent $\alpha_{1}$ d\'esignant, comme au d\'ebut du paragraphe 3.3, le plus petit \'el\'ement du spectre de $(G^{o},G_{0}^{o})$ et si la multiplicit\'e de $\alpha_{1}$ dans ce spectre est \'egale \`a $1$
(voir [D, proposition 7.0.2]). La premi\`ere condition signifie que $\phi^{*}(1)$ est {\em minimal} : assez curieusement,
et contrairement \`a ce qui se passe dans le cas local, ce n'est pas toujours le cas (voir la remarque ci-dessous). 
Si $f$ est quasi-homog\`ene, on obtient cependant toujours une forme homog\`ene $\zeta$ de cette fa\c{c}on.\\

\noindent $\bullet$ La condition d'injectivit\'e $(IC)$ de [HeMa, theorem 4.5] se lit, \'etant entendu que la forme $\zeta$ choisie est donn\'ee par la classe de la forme volume
$du:=du_{1}\wedge\cdots\wedge du_{n}$,
\begin{center} {\em les classes de $u_{1}du,\cdots ,u_{n}du$ dans $\Omega^{n}(\cit^{n})/df\wedge\Omega^{n-1}(\cit^{n})$ sont lin\'eairement ind\'ependantes}.
\end{center}
Cette condition d'injectivit\'e est toujours garantie si l'on fait l'hypoth\`ese 
$$\phi^{*}(u_{j})<1-\phi (u_{i})$$
\noindent pour tout $i,j\in \{1,\cdots ,n\}$ : ceci se voit comme au lemme 4.1.1. Cette condition impose $\phi (u_{i})<1$ pour tout $i$. 

\begin{remark} Dans le cas local, {\em i.e} si $f$ est une fonction holomorphe sur un voisinage ouvert de l'origine, \`a point critique isol\'e en $0$, on peut d\'efinir de mani\`ere analogue une filtration de Newton $\phi$ sur 
$\cit \{u_{1},\cdots ,u_{n}\}$ (voir [K]). On a alors $\phi (ab)\geq \phi (a)+\phi (b)$ et donc 
$$\phi^{*}(g)\geq \phi^{*}(1)+\phi (g)$$
\noindent pour tout $g\in \cit \{u_{1},\cdots ,u_{n}\}$ : ceci montre que $\phi^{*}(1)$ est minimal. Si $g\in\cit [u_{1},\cdots ,u_{n}]$ on a seulement
$$\phi^{*}(g)\leq \phi^{*}(1)+\phi (g).$$
\end{remark}

\noindent Une fois ces deux difficult\'es surmont\'ees, on peut recopier ce qui a \'et\'e fait en section 4.1 et obtenir  

\begin{theorem}  Soit $f$ une application polynomiale non d\'eg\'en\'er\'ee et commode. On suppose que
$\phi^{*}(u_{j})<1-\phi (u_{i})$ (et donc $\phi (u_{i})<1$) pour tout $i,j\in \{1,\cdots ,n\}$. Alors
la trTLEP(n)-structure fournie, selon le proc\'ed\'e du paragraphe 3, par le r\'eseau de Brieskorn de la d\'eformation 
$$F(u,x)=f(u)+\sum_{j=1}^{n}x_{j}u_{j}$$
admet une d\'eformation universelle.
\end{theorem}

\begin{remark} $1.$ L'hypoth\`ese sur les $\phi^{*}(u_{j})$ n'est faite que pour garantir la condition (IC) : elle est souvent trop forte (c'est le cas dans l'exemple ci-dessous). Par contre, et \`a cause de la proposition 4.2.1, notre m\'ethode requiert toujours $\phi (u_{i})<1$.\\
$2.$ Le choix de la d\'eformation consid\'er\'ee est motiv\'e par la condition (GC) : bien entendu, il n'est pas toujours optimal. Par exemple, si l'on consid\`ere $f(u_{1},u_{2})=u_{1}^{3}+u_{2}^{2}$ il suffit de prendre $F(u_{1},u_{2},x_{1})=u_{1}^{3}+u_{2}^{2}+x_{1}u_{1}$.
\end{remark}

\noindent Cette mani\`ere de voir est particuli\`erement utile si l'on consid\`ere les polyn\^omes quasi-homog\`enes, polyn\^omes pour lesquels on a toujours $R_{0}=0$ sur $E_{0}$ : dans la condition (GC), 
l'hypoth\`ese portant sur $R_{0}$ est par vide. Le th\'eor\`eme 4.2.5 est par cons\'equent 
id\'ealement adapt\'e \`a cette situation.

\subsubsection{Un exemple}

Soit $f:\cit^{2}\rightarrow\cit$, $f(u_{1},u_{2})=u_{1}^{5}+u_{2}^{5}$. On a $\mu =16$ et les mon\^omes 
$$1,u_{1},u_{2}, u_{1}^{2},u_{1}u_{2},u_{2}^{2},u_{1}^{3},u_{1}^{2}u_{2}, u_{1}u_{2}^{2},u_{2}^{3}, u_{1}^{3}u_{2}, u_{1}^{2}u_{2}^{2},u_{1}u_{2}^{3},u_{1}^{3}u_{2}^{2}, u_{1}^{2}u_{2}^{3}, u_{1}^{3}u_{2}^{3}$$
forment une base de $\cit [u_{1},u_{2}]/(\partial f/\partial u_{1},\partial f/\partial u_{2})$.
Un d\'eploiement universel de $f$ peut s'\'ecrire
$$\tilde{F}(u_{1},u_{2},x)=f(u_{1},u_{2})+\sum_{i=1}^{16}x_{i}g_{i}(u)$$
o\`u les $g_{i}$ parcourent l'ensemble des mon\^omes ci-dessus. Cet exemple est particuli\`erement int\'eressant parce qu'il s'agit d'un des cas le plus simple o\`u $\tilde{F}$ a des points critiques qui disparaissent \`a l'infini quand les param\`etres $x_{i}$ tendent vers $0$ et pour lequel par cons\'equent notre m\'ethode est particuli\`erement efficace.\\

\noindent Avec les notations ci-dessus, on a
$$\phi ^{*}(u_{1}^{i}u_{2}^{j})=\frac{i+j+2}{5}$$
\noindent pour tout $i\in\nit$ et tout $j\in\nit$. En particulier, $\phi^{*}(1)=2/5$ est minimal et cette valeur n'est atteinte que pour le mon\^ome $1$ : la classe de la forme $du_{1}\wedge du_{2}$ est donc primitive et homog\`ene. On a aussi $\phi (u_{1})=1/5 <1$ et
$\phi (u_{2})=1/5<1$.\\

\noindent On consid\`ere donc
$$F(u_{1},u_{2},x_{1},x_{2})=f(u_{1},u_{2})+x_{1}u_{1}+x_{2}u_{2}.$$
Il est clair que les classes de $u_{1}$ et $u_{2}$ dans $\cit [u_{1},u_{2}]/(\partial_{u_{1}}f,\partial_{u_{2}}f)$ sont ind\'ependantes. La condition $(IC)$ est par cons\'equent v\'erifi\'ee. La condition $(GC)$ r\'esulte du choix des mon\^omes $u_{1}$ et $u_{2}$.
La seule chose \`a faire est de r\'esoudre le probl\`eme de Birkhoff pour le r\'eseau de Brieskorn $G_{0}$ associ\'e \`a $F$. Il nous faut donc trouver une base $\epsilon$ qui v\'erifie les conclusions du corollaire 3.1.3 et qui se comporte bien vis-\`a-vis de la dualit\'e. Tout ceci peut se faire \`a la main.\\
 
\noindent Notons $\epsilon_{i,j}=[u_{1}^{i}u_{2}^{j}du_{1}\wedge du_{2}]$ o\`u $[\ ]$ d\'esigne la classe de l'\'el\'ement consid\'er\'e dans le r\'eseau de Brieskorn $G_{0}$. Comme $\phi (u_{1})<1$ et $\phi (u_{2})<1$, la proposition 2.2.1 nous fournit les r\'esultats suivants :\\

$\bullet$ $G_{0}$ est libre, de rang $16$, sur $\cit [x_{1},x_{2},\theta ]$,\\

$\bullet$ $\epsilon =(\epsilon_{i,j})_{0\leq i\leq 3, 0\leq j\leq 3}$ est une base de $G_{0}$ sur $\cit [x_{1},x_{2},\theta ]$.\\

\noindent Nous allons v\'erifier tout d'abord que cette base est bien une solution au probl\`eme de Birkhoff au sens du corollaire 3.1.3.
On a
$$Fu_{1}^{i}u_{2}^{j}=\frac{4}{5}x_{1}u_{1}^{i+1}u_{2}^{j}+\frac{4}{5}x_{2}u_{1}^{i}u_{2}^{j+1}+dF\wedge [-\frac{1}{5}u_{1}^{i}u_{2}^{j+1}du_{1}+\frac{1}{5}u_{1}^{i+1}u_{2}^{j}du_{2}]$$
d'o\`u 
$$\theta^{2}\nabla_{\partial_{\theta}}\epsilon_{i,j}=
\frac{4}{5}x_{1}\epsilon_{i+1,j}+\frac{4}{5}x_{2}\epsilon_{i,j+1}
+\frac{i+j+2}{5}\theta \epsilon_{i,j}$$
dans $G_{0}$ pour $0\leq i\leq 2$ et $0\leq j\leq 2$.\\
De plus,
$$Fu_{1}^{3}u_{2}^{j}=-\frac{4}{25}x_{1}^{2}u_{2}^{j}+\frac{4}{5}x_{2}u_{1}^{3}u_{2}^{j+1}+
dF\wedge [-\frac{1}{5}u_{1}^{3}u_{2}^{j+1}du_{1}+\frac{1}{5}(u_{1}^{4}u_{2}^{j}+\frac{4}{5}x_{1}u_{2}^{j})du_{2}],$$

$$Fu_{1}^{i}u_{2}^{3}=-\frac{4}{25}x_{2}^{2}u_{1}^{i}+\frac{4}{5}x_{1}u_{1}^{i+1}u_{2}^{3}+
dF\wedge [-\frac{1}{5}(u_{1}^{i}u_{2}^{4}+\frac{4}{5}x_{2}u_{1}^{i})du_{1}+
\frac{1}{5}u_{1}^{i+1}u_{2}^{3}du_{2}]$$
et 
$$Fu_{1}^{3}u_{2}^{3}=-\frac{4}{25}x_{1}^{2}u_{2}^{3}-\frac{4}{25}x_{2}^{2}u_{1}^{3}+
dF\wedge [(-\frac{1}{5}u_{1}^{3}u_{2}^{4}-\frac{4}{25}x_{2}u_{1}^{3})du_{1}+
(\frac{1}{5}u_{1}^{4}u_{2}^{3}+\frac{4}{25}x_{1}u_{2}^{3})du_{2}].$$
On en d\'eduit
$$\theta^{2}\nabla_{\partial_{\theta}}\epsilon_{3,j}=
-\frac{4}{25}x_{1}^{2}\epsilon_{0,j}+\frac{4}{5}x_{2}\epsilon_{3,j+1}
+\frac{5+j}{5}\theta \epsilon_{3,j}$$
pour $0\leq j\leq 2$,
$$\theta^{2}\nabla_{\partial_{\theta}}\epsilon_{i,3}=
-\frac{4}{25}x_{2}^{2}\epsilon_{i,0}+\frac{4}{5}x_{1}\epsilon_{i+1,3}
+\frac{5+i}{5}\theta \epsilon_{i,3}$$
pour $0\leq i\leq 2$
et 
$$\theta^{2}\nabla_{\partial_{\theta}}\epsilon_{3,3}=
-\frac{4}{25}x_{1}^{2}\epsilon_{0,3}-\frac{4}{25}x_{2}^{2}\epsilon_{3,0}
+\frac{8}{5}\theta \epsilon_{3,3}.$$
Ceci montre que la matrice de $\theta^{2}\nabla_{\partial_{\theta}}$ dans la base $\epsilon$ s'\'ecrit bien sous la forme voulue : les matrices $B_{0}(x)$ et $B_{\infty}$ se d\'eduisent des formules pr\'ec\'edentes. Remarquons que $B_{\infty}$ est constante, diagonale, et que ses valeurs propres sont les rationnels $\phi ^{*}(u_{1}^{i}u_{2}^{j})$, $0\leq i\leq 3$, $0\leq j\leq 3$.\\

\noindent Des calculs analogues aux pr\'ec\'edents (il suffit d'\'evaluer la multiplication par $-u_{i}$, $i=1,2$, dans la base $\epsilon$) montrent que la matrice de $\nabla_{\partial_{x_{i}}}$, $i=1,2$, dans la base $\epsilon$ 
s'\'ecrit $C^{(i)}(x)\theta^{-1}$ avec 
$$C^{(1)}(x)(\epsilon_{i,j})=-\epsilon_{i+1,j}$$
si $i\leq 2$,
$$C^{(1)}(x)(\epsilon_{3,j})=+\frac{x_{1}}{5}\epsilon_{0,j},$$
$$C^{(2)}(x)(\epsilon_{i,j})=-\epsilon_{i,j+1}$$
si $j\leq 2$,
$$C^{(2)}(x)(\epsilon_{i,3})=\frac{x_{2}}{5}\epsilon_{i,0}.$$

\noindent En r\'esum\'e, la matrice de la connexion $\nabla$ dans la base $\epsilon$ s'\'ecrit
$$(\frac{B_{0}(x)}{\theta}+B_{\infty})\frac{d\theta}{\theta}+\frac{C^{(1)}(x)}{\theta}dx_{1}+
\frac{C^{(2)}(x)}{\theta}dx_{2}.$$

\noindent Terminons par le comportement de la base $\epsilon$ vis-\`a-vis de la dualit\'e : on pose
$$S(\epsilon_{i,j},\epsilon_{k,l})=\tau^{-n}$$
si $k=3-i$ et $l=3-j$,
$$S(\epsilon_{i,j},\epsilon_{k,l})=0$$
sinon. Ceci fournit une forme sequilin\'eaire sur $G$ (voir le paragraphe 3.3.3). Les calculs ci-dessus montrent que les matrices $B_{\infty}$,$B_{0}(x)$, $C^{(1)}(x)$ et $C^{(2)}(x)$ sont bien sym\'etriques au sens o\`u 
$$B_{\infty}+B_{\infty}^{*}=nI,\  B_{0}(x)^{*}=B_{0}(x),\ C^{(1)}(x)^{*}=C^{(1)}(x),\ C^{(2)}(x)^{*}=C^{(2)}(x),$$ 
\noindent $B^{*}$ d\'esignant l'adjoint de $B$ relativement \`a $S$. Par cons\'equent, $S$ est bien compatible aux connexions, de poids $n$ relativement au r\'eseau $G_{0}$. La solution au probl\`eme de Birkhoff trouv\'ee ici ne d\'epend pas des param\`etres $x_{1}$ et $x_{2}$ : ceci est un cas tout \`a fait exceptionnel.

\section{$\mu$-constant : remarques sur le cas g\'en\'eral}

La d\'eformation envisag\'ee dans les paragraphes 2 et 3 de ces notes peut sembler restrictive, m\^eme si elle suffit 
si l'on privil\'egie la construction de structures de Frobenius.
Il serait tr\`es int\'eressant d'\'etudier le syst\`eme de Gauss-Manin et le r\'eseau de Brieskorn 
associ\'es \`a une famille g\'en\'erale de 
polyn\^omes $M$-tame (et donc plus n\'ecessairement non d\'eg\'en\'er\'es) \`a nombre de Milnor global constant et de g\'en\'eraliser par exemple les 
r\'esultats de la proposition 2.3.1. 
Dans le cas local c'est une affaire entendue (voir par exemple [Va] et [Sai]). Cependant, 
m\^eme dans le cas commode et non d\'eg\'en\'er\'e, les arguments employ\'es dans ces notes sont mis en d\'efaut 
(le lemme de division n'apporte rien dans les cas 'limites'). Se pose en particulier le probl\`eme de l'identification de la $V$-filtration.
Il faut donc imaginer autre chose.\\

{\bf Coh\'erence.} 
La premi\`ere question qui se pose est la coh\'erence du r\'eseau de Brieskorn $G_{0}$ d'une telle d\'eformation : ceci devrait \^etre fait 
 en utilisant par exemple les arguments expos\'es dans [BG], 
en particulier le th\'eor\`eme de Kiehl-Verdier. Il faut cependant faire attention au ph\'enom\`ene \'eventuel des 
``plis \'evanescents (\`a l'infini)'' et utiliser les r\'esultats de L\^e-Ramanujam [LR] qui 
s'\'etendent sans difficult\'es, avec les m\^emes preuves, aux polyn\^omes $M$-tame 
(voir [VZ]). On obtient ainsi un 
$\cal{O}_{M}[t]$-module de type fini $M_{0}$ et, par transformation de Fourier un $\cal{O}_{M}[\theta ]$-module libre de 
type fini $G_{0}$ en \'evitant le 
proc\'ed\'e d'alg\'ebrisation de Malgrange utilis\'e dans [DoSa].\\

{\bf Fibr\'es et sous-fibr\'es.}
La deuxi\`eme question, qui s'impose naturellement si l'on veut \'etendre les r\'esultats obtenus ci-dessus 
(en particulier ceux de la proposition 2.3.1), 
 est la suivante : on doit s'attendre (et c'est l\`a un point clef qui devrait pouvoir se montrer en remarquant que $G$ est 
l'image r\'eciproque d'un
$\cit [\tau ]<\partial_{\tau}>$-module holonome) \`a ce que les 
$$H_{\alpha}:=V_{\alpha}G/V_{<\alpha}G$$
soient des $\cal{O}_{M}$-modules localement libres 
 et
d\'efinir de la m\^eme mani\`ere qu'en 2.3 (voir plus particuli\`erement le corollaire 2.3.2), 
\`a l'aide de la $V$-filtration le long de $\tau :=\theta^{-1}=0$, une filtration de chaque $H_{\alpha}$ par des 
$\cal{O}_{M}$-modules 
$$F_{p}H_{\alpha}:=\frac{V_{\alpha}G\cap\tau^{p}G_{0}+V_{<\alpha }G}{V_{<\alpha }G}$$ 
(la filtration de Hodge). Le point est maintenant le suivant : 
les $F_{p}H_{\alpha}$ sont-ils des sous-fibr\'es   
(peut on esp\'erer obtenir sur $M$ une variation de structures de Hodge mixtes au sens de Steenbrink-Zucker)? La r\'eponse est oui dans le cas local
(c'est ce que montre Varchenko dans [Va]). Signalons que cette question ({\em les $F_{p}$ sont-ils des sous-fibr\'es?})
est ind\'ependante de la r\'esolution du probl\`eme de Birkhoff (et donc de la construction de structures de Frobenius). Remarquons aussi, pour faire 
bonne mesure, que la condition de transverslit\'e de 
Griffiths se lit {\em a priori}
$$\nabla_{\partial_{x}} F_{p}H_{\alpha}\subset F_{p+1}H_{\alpha}$$
(\`a comparer avec la remarque 2.3.2).\\

{\bf Le probl\`eme de Birkhoff.} 
Si les $H_{\alpha}$ sont $\cal{O}_{M}$-localement libres, le probl\`eme de Birkhoff pour $G_{0}$ admet une solution ais\'ee (par voie 
{\em analytique} cependant), une fois que l'on sait
le r\'esoudre pour la valeur $0$ des param\`etres (il est remarquable ici de constater que l'on peut oublier la filtration $F_{\bullet}$ sur $H_{\alpha}$,
il suffit juste de la consid\'erer pour la valeur nulle des param\`etres) : on utilise
la d\'ecomposition 
$$H_{\alpha}^{0}=\oplus_{p\in\zit}F_{p}H_{\alpha}^{0}\cap U^{p}H_{\alpha}^{0}$$
($H_{\alpha}^{0}$ d\'esigne comme ci-dessus la restriction de $H_{\alpha}$ \`a la valeur nulle des param\`etres)
pour obtenir, par transport parall\`ele, la d\'ecomposition
$$H_{\alpha}=\oplus_{p\in\zit}G_{\alpha ,p}$$
(les  $G_{\alpha ,p}$ sont des $\cal{O}_{M}$-modules localement libres, rappelons que $H_{\alpha}$ est suppos\'e \^etre $\cal{O}_{M}$-localement libre).
Cette d\'ecomposition permet de d\'efinir une filtration d\'ecroissante $U^{\bullet}$
de $H_{\alpha}$, stable par $\tau\nabla_{\partial_{\tau}}$ et $\nabla_{\partial_{x_{i}}}$ (de la m\^eme mani\`ere qu' au paragraphe 3) en posant
$$U^{p}H_{\alpha}=\oplus_{q\geq p}G_{\alpha ,q}.$$
A son tour, par le proc\'ed\'e maintenant habituel, cette filtration d\'ecroissante donne un r\'eseau $G_{\infty}$, logarithmique le long de 
$\tau =0$ (ou encore $G_{\infty}$ v\'erifie les conditions {\bf b.} et {\bf c.} de la proposition 3.1.1). Ce r\'eseau permet de d\'efinir
un fibr\'e trivial sur $\ppit^{1}\times M$ (quitte \`a r\'etr\'ecir $M$), parce que la trivialit\'e est une propri\'et\'e ouverte, muni d'une connexion
\`a p\^oles logarithmiques le long de $\tau =0$ et de rang de Poincar\'e inf\'erieur ou \'egal \`a $1$ en $\tau =\infty$. Ceci donne une solution 
au probl\`eme de Birkhoff pour $G_{0}$. Notons cependant que, 
contrairement au cas trait\'e dans les paragraphes 2. et 3. de 
ces notes, la restriction de cette solution n'est pas \'egale \`a la solution de la restriction donn\'ee par [DoSa, appendix B], 
parce qu'en g\'en\'eral $\nabla_{\partial_{x}} F_{p}H_{\alpha}\subset F_{p+1}H_{\alpha}$ : si tel \'etait le cas, on aurait, 
avec les notations pr\'ec\'edentes,
$F_{p}H_{\alpha}=\oplus_{q\leq p}G_{\alpha ,q}$ et donc 
$$\nabla_{\partial_{x}} F_{p}H_{\alpha}\subset F_{p}H_{\alpha}.$$
Remarquons que cette derni\`ere inclusion est vraie si et seulement si, pour tout $i$ et tout $\omega\in \Omega^{n}(U)$ dont la classe appartient \`a 
$V_{\alpha}G\cap G_{0}$, la classe de 
$\frac{\partial F}{\partial_{x_{i}}}\omega$ dans $G$ appartient en fait \`a
$$V_{<\alpha +1}G\cap G_{0}+V_{\alpha +1}G\cap\theta G_{0}.$$
Quoi qu'il en soit, il n'est pas ais\'e, except\'e dans des cas simples comme celui consid\'er\'e dans les paragraphes 2 et 3 de ces notes, 
de le v\'erifier, ne serait-ce que parce que l'on ne dispose pas de description explicite de la $V$-filtration.  

Nous sommes maintenant en mesure d'affiner [Mal, th\'eor\`eme 2.2]. Soit en effet $(\omega_{1},\cdots ,\omega_{\mu})$ une 
solution du probl\`eme de 
Birkhoff pour la valeur nulle des param\`etres donn\'ee par [DoSa, appendix B], avec, pour tout $i$,
 $$\omega_{i}\in V_{\alpha_{i}}G^{o}\cap G_{0}^{o}$$
$G^{o}$ ({\em resp.} $G_{0}^{o}$) d\'esignant le syst\`eme de Gauss-Manin ({\em resp.} le r\'eseau de Brieskorn) absolu et $V_{\bullet}$ 
la filtration de Malgrange-Kashiwara.
 Si $(\tilde{\omega}_{1},\cdots ,\tilde{\omega}_{\mu})$ d\'esigne l'extension de cette solution donn\'ee par le r\'esultat de Malgrange, alors la restriction
$\tilde{\omega}^{x}_{i}$ de la section globale $\tilde{\omega}_{i}$ \`a la valeur $x$ des param\`etres appartient \`a  
$$V_{\alpha_{i}}G^{x}\cap G_{0}^{x},$$
 $G^{x}$ ({\em resp.} $G_{0}^{x}$) d\'esignant la restriction \`a $x$ du syst\`eme de Gauss-Manin $G$ 
({\em resp.} du r\'eseau de Brieskorn $G_{0}$) 
associ\'e \`a la d\'eformation \`a $\mu$-constant consid\'er\'ee (il s'agit en fait des transform\'es de Fourier du syst\`eme de Gauss-Manin et du 
r\'eseau de Brieskorn habituels, voir le d\'ebut du paragraphe 2). 
Ce r\'esultat doit \^etre rapproch\'e de [Sai, (3.6.3)]; 
on peut d'ailleurs penser que l'on obtient ainsi une carat\'erisation de la strate \`a $\mu$ (global)-constant.

Signalons pour finir que le comportement par dualit\'e ne pose pas 
de probl\`emes : si l'on part d'une $S$-solution (au sens de [DoSa, appendix B]) on obtient par ce proc\'ed\'e une $S$-solution parce que la 
dualit\'e est compatible avec les connexions.

\begin{flushright}
Antoine DOUAI \\
Universit\'e de Nice \\
UMR 6621 du CNRS \\
Laboratoire J. A. Dieudonn\'e \\
06 108 Nice cedex 2 \\
douai@math.unice.fr \\
http://math.unice.fr/$\sim$ douai 
\end{flushright}


\begin{thebibliography}{999} 
\bibitem[BG]{BG} Buchweitz (R.A), Greuel (G.M).--{\em  The Milnor number and deformations of complex 
curve singularities}, Inv. Math., t {\bf 58}, 1980.
\bibitem[D]{D} Douai (A).--{\em  Contributions \`a l'\'etude des syst\`emes de Gauss-Manin alg\'ebriques}, pr\'epublication 695 de l'universit\'e de Nice, 2004. Disponible sur le site personnel de l'auteur.
\bibitem[DoSa]{DoSa} Douai (A), Sabbah (C).-- {\em Gauss-Manin systems, Brieskorn lattices and Frobenius structures I}, Ann. Inst. Fourier, 53-4 (2003), 
1055-1116.
\bibitem[DoSa2]{DoSa2} Douai (A), Sabbah (C).-- {\em Gauss-Manin systems, Brieskorn lattices and Frobenius structures II}, in {\em Frobenius Manifolds},
C. Hertling and M. Marcolli (Eds.), Aspects of Mathematics E 36.
\bibitem[He]{He} Hertling (C).--{\em Frobenius manifolds and moduli spaces for singularities}, Cambridge tracts in Math., {\bf 151}, 2002.
\bibitem[HeMa]{HeMa} Hertling (C), Manin (Y).--{\em Unfoldings of meromorphic connections and a construction of Frobenius manifolds}, in {\em Frobenius Manifolds}, C. Hertling and M. Marcolli (Eds), Aspects of Mathematics E 36.
\bibitem[K]{K} Kouchnirenko (A.G).--{\em Poly\`edres de Newton et nombres de Milnor}, Inv. Math., t. {\bf 32}, 1976, p. 1-31.
\bibitem[LR]{LR} L\^e (D.T), Ramanujam (C.P).--{\em The invariance of Milnor's number implies 
the invariance of the topological type}, American journal of math., t. {\bf 98}, 1976.
\bibitem[Mal]{Mal} Malgrange (B).--{\em D\'eformations de syst\`emes diff\'erentiels et microdiff\'erentiels}, dans {\em S\'eminaire E.N.S Math\'ematique
 et Physique}, t. {\bf 6}, p. 351-379.
\bibitem[Ph]{Ph} Pham (F).-- {\em Singularit\'es des syst\`emes diff\'erentiels de Gauss-Manin}, Progress in Math., t. {\bf 2}, Birkhauser, Boston, 1980.
\bibitem[Sab1]{Sab1} Sabbah (C).-- {\em Hypergeometric periods for a tame  polynomial}, C.R Acad. Sci. Paris S\'e. I Math. , t. {\bf 328}, 1999, 
p. 603-608 et preprint math. AG/9805077.
\bibitem[Sab2]{Sab2} Sabbah (C).-- {\em Monodromy at infinity and Fourier transform}, Publ. RIMS, Kyoto Univ., t. {\bf 33}, 1998, p. 643-685. 
\bibitem[Sab3]{Sab3} Sabbah (C).--{\em D\'eformations isomonodromiques et vari\'et\'es de Frobenius}, Savoirs Actuels, CNRS Editions, Paris, 2002.
\bibitem[Sai]{Sai} Saito (M).--{\em  Period mapping via Brieskorn modules}, Bull. Soc. Math. France, t. {\bf 119}, 1991, p. 141-171.
\bibitem[Sai2]{Sai2} Saito (M).--{\em  On the structure of Brieskorn lattices}, Ann. Inst. Fourier, t. {\bf 39}, 1989, p. 27-72.
\bibitem[Va]{Va} Varchenko (A).--{\em The complex exponent of a singularity does not change along strata $\mu =$ const}, Functionnal Anal. Appl., 
t. {\bf 16},1982, p. 469-512.
\bibitem[VZ]{VZ} Vui (H.H), Zaharia (A).--{\em Families of polynomials with total Milnor number}, Math. Ann., t. {\bf 304}, 1996, p. 481-488.


\end{thebibliography}
\end{document}